\documentclass[12pt]{article}
\usepackage{amssymb,amsmath,amsopn,amsfonts}
\usepackage[dvips]{color}
\usepackage{colordvi,multicol}
\usepackage{epsf}
\newfam\msbfam
\font\tenmsb=msbm10 \textfont\msbfam=\tenmsb \font\sevenmsb=msbm7
\scriptfont\msbfam=\sevenmsb \font\fivemsb=msbm5
\scriptscriptfont\msbfam=\fivemsb

\def\th#1{\vspace{1mm}\noindent{\bf #1}\quad}

\def\proof{\vspace{1mm}\noindent{\it Proof}\quad}

\numberwithin{equation}{section}

\def\bc{\begin{center}}
\def\ec{\end{center}}
\def\no{\noindent}
\def\hang{\hangindent\parindent}
\def\textindent#1{\indent\llap{\qquad #1\ \ \enspace}\ignorespaces}
\def\ref{\par\hang\textindent}

\begin{document}

\title{ {\bf Local existence and non-explosion of solutions for  stochastic fractional partial differential equations driven by multiplicative noise
\thanks{Research supported in part by China Postdoctoral Science Foundation funded project(2012M520153) and the DFG through IRTG 1132 and CRC 701 }\\} }
\author{ {\bf Michael R\"{o}ckner}$^{\mbox{c}}$, {\bf Rongchan Zhu}$^{\mbox{b,c},}$, {\bf Xiangchan Zhu}$^{\mbox{a,c},}$
\thanks{E-mail address: roeckner@math.uni-bielefeld.de(M. R\"{o}ckner),
zhurongchan@126.com(R. C. Zhu), zhuxiangchan@126.com(X. C. Zhu)}\\\\
$^{\mbox{a}}$School of Sciences, Beijing Jiaotong University, Beijing 100044, China, \\
$^{\mbox{b}}$Department of Mathematics, Beijing Institute of Technology, Beijing 100081,  China
 \\
$^{\mbox{c}}$Department of Mathematics, University of Bielefeld, D-33615 Bielefeld, Germany}
\maketitle

\maketitle

\noindent {\bf Abstract}

In this paper we prove  the local existence and uniqueness of solutions for a class of stochastic fractional partial differential equations driven by multiplicative noise.
We also establish that for this class of equations adding linear multiplicative noise provides a
regularizing effect: the solutions will not blow up with high probability if the initial data is sufficiently
small, or if the noise coefficient is sufficiently large. As applications our main results are applied to various types of SPDE such as stochastic reaction-diffusion equations, stochastic fractional Burgers equation, stochastic fractional Navier-Stokes equation, stochastic quasi-geostrophic equations and stochastic surface growth PDE.

\vspace{1mm}
\no{\footnotesize{\bf 2000 Mathematics Subject Classification AMS}:\hspace{2mm} 60H15, 60H30, 35R60}
 \vspace{2mm}

\no{\footnotesize{\bf Keywords}:\hspace{2mm}   stochastic fractional partial differential equation, local existence and uniqueness, blow up, Navier-Stokes equation, fractional Burgers equation, quasi-geostrophic equation, surface growth models }

\section{Introduction}
Consider the following  stochastic equation with fractional dissipation in a smooth bounded domain $O\subset \mathbb{R}^d$:
$$\frac{\partial \theta(t,
\xi)}{\partial t}+f(\theta)(t,\xi)+A^\alpha \theta(t,\xi)=(G(\theta)\eta)(t,\xi),\eqno(1.1)$$
with  initial condition $$\theta(0,\xi)=\theta_0(\xi), \eqno(1.2)$$
where $\theta(t,\xi)$ is a (possibly vector valued) function of $\xi\in O$ and $t\geq0$, $0<\alpha\leq1$ are real numbers.
Here $A$ is an unbounded positive definite  self-adjoint operator, $f(\theta)$ is a nonlinear term and $\eta(t,\xi)$ is a Gaussian random field, white noise in time, subject to the restrictions imposed below.

This model includes a large number of nonlinear models in fluid dynamics such as  the fractional Burgers equation (see [KNS08]),
  the quasi-geostrophic equation (see [CV10], [RZZ12] and the reference therein), the Navier-Stokes equation (see e.g.[GF95]),  surface growth models (see e.g.[BFR09]) and reaction-diffusion equations.

 The main result in [LR13] establishes the local existence and uniqueness of solutions to (1.1) driven by additive noise if the coefficients satisfy some local monotonicity and generalized coercivity conditions. However, the supercritical case (i.e. $\alpha<\frac{1}{2}$) of the stochastic quasi-geostrophic equation and the stochastic fractional Burgers equation are not within the framework in [LR13]. To the best of our knowledge there are no results in the literature that cover these cases.

In the first part of this paper we establish local existence and uniqueness of solutions in $C([0,\infty);H^{s_0})$ (for the definition of $H^{s_0}$ see below) for equation (1.1) (see Theorem 3.2). Here we emphasize that our results in particular cover the supercritical case (i.e. $\alpha<\frac{1}{2}$) of the quasi-geostrophic equation and the fractional Burgers equation. In particular,  using our result for the critical case (i.e. $\alpha=\frac{1}{2}$) of the quasi-geostrophic equation driven by linear multiplicative noise we obtain existence and uniqueness of global solutions (see Remark 5.7).
The technical key point to achieve this is to identify the appropriate coercivity condition (b.1) and the local monotonicity condition (b.3) which are formulated to hold in two different Hilbert spaces for the existence and uniqueness respectively. Moreover, due to the dissipation term of the equation, by a boot-strapping argument, we can  deduce the solution has enough regularity to control the nonlinear term (see (b.2)). We also emphasize that our conditions are satisfied by all the examples mentioned above (see Section 5).

The first difficulty arising is to prove that the solution is continuous with respect to $H^{s_0}$-norm when the initial value is in $H^{s_0}$. In [GV12] the authors prove  local existence and uniqueness of solutions for the stochastic Euler equation. In that paper the result is first proved for smooth initial values and then the authors deal with more general initial values by approximation. However, this method cannot be applied to the fractional Burgers equation  since we cannot use the commutator estimate (Lemma 5.3) in $H^s$ to control the nonlinear term of the fractional Burgers equation. Instead, we use a boot-strapping argument to deduce that the solution is in $C((0,T];H^{s_0})$ and then we prove the continuity in $t=0$. The second obstacles lies in the lack of the Yamada-Watanabe theorem for the local solution. Here we view the local solution as a global solution for another equation (see (3.18)) and  hence we can use the Yamada-Watanabe theorem from [Ku07] for equation (3.18) to obtain the existence and uniqueness of a (probabilistically) strong solution.

Recently, there has been a lot of work done on how results on
  partial differential equations (PDE) change due to random perturbations (see [FGP10, Fl10, DFPR12, GV12] and references therein). A very interesting case of regularization by noise is
the case when the PDE is not well posed, but the SPDE is well posed. In the deterministic case, when $\alpha<1/2$ the global existence and uniqueness of smooth solutions for the quasi-geostrophic equation remains an open problem. Also the uniqueness of weak solutions for the surface growth model is still an open problem. For the fractional Burgers equation in the supercrtical case
the solution will  blow up in $H^s$ in finite time when $s>\frac{3}{2}-2\alpha$ (see [KNS08, Theorem 1.4]). In [GV12] the authors obtain that adding  linear multiplicative noise to the 3D Euler equation provides a regularizing effect. Inspired by this we establish such kind of result for (1.1). More precisely:

In the second part of this paper  we consider a special case of a linear multiplicative noise $\beta \theta dW$ with $W$ a real-valued Brownian motion, $\beta\in\mathbb{R}$, and prove that  for any $R\geq1$: whenever $\|\theta_0\|_{H^{s_0}}\leq \kappa(R,\beta^2),$

$$P(\theta \textrm{ is global })\geq1-R^{-1/4}, $$and
$$P(\|\theta(t)\|_{H^{s_0}}\rightarrow0 \textrm{ as } t\rightarrow\infty)\geq1-R^{-1/8}, $$
where $\kappa=\kappa(R,\beta^2)$ is some explicit functions such that $$\lim_{\beta^2\rightarrow\infty}\kappa(R,\beta^2)=\infty,$$
for every fixed $R\geq 1$ (see Theorem 4.1). This can be viewed as a kind of global existence result in the large noise asymptotics.
We transform (1.1) to a PDE with a random coefficient producing a damping term. We can exploit this random damping as done in [GV12]. This damping term is strong enough to force the solution to go to $0$ as $t\rightarrow\infty$.

This paper is organized as follows. In Section 2 we recall some preliminaries. In Section 3, we  prove  the existence and uniqueness  of local solutions to (1.1). The regularization result obtained through perturbing by noise is given in Section 4.  In Section 5, we apply the results from Sections 3 and 4 to the concrete fractional SPDE mentioned at the beginning of this introduction.

\section{Notations and Preliminaries}

 Let $H$ be a separable Hilbert space  and $|\cdot|$, $\langle \cdot,\cdot\rangle$ denote the norm and inner product in $H$ respectively. Let $A:D(A)\subset H\rightarrow H$ be a positive definite  self-adjoint operator such that $A^{-1}$ is compact on $H$. From this we infer the existence of a complete orthonormal basis $\{e_k\}_{k\geq0}$ for $H$ of eigenfunctions of $A$ such that the associated sequence of eigenvalues $\{\lambda_k\}$ form an increasing unbounded sequence.

 Using the basis $\{e_k\}$ we may also define the fractional powers of $A$. Given $s>0$ define
 $$H^s:=D(A^{\frac{s}{2}})=\{f\in H:\sum_k\lambda_k^{s}|\langle f,e_k\rangle|^2<\infty\},$$
 and $$A^{\frac{s}{2}} f:=\sum_k\lambda_k^{\frac{s}{2}}\langle f,e_k\rangle e_k, \quad f\in D(A^s).$$
 For $s<0$ define $H^s$ to be the dual of $H^{-s}$.
 Set $\Lambda=A^{1/2}$. For $s>0$, define the associated Hilbert norm
$$\|f\|_{H^s}^2:=|\Lambda^s f|^2=\sum_k\lambda_k^s|\langle f,e_k\rangle|^2.$$
Then we can easily verify that $H^{s_1}\subset H^{s_2}$ is a compact embedding if $s_1>s_2$.
 Let $\Phi(H^r;H^u)$ denote all measurable mappings from $H^r$ to $H^u$ for $r,u\in\mathbb{R}$ and $H_n=\textrm{span}\{e_j\}_{j\leq n}$.

First we prove the following lemma for later use.

\th{Lemma 2.1} $B\subset C((0,T];H^s)$ is relatively compact with respect to the topology induced by the complete metric $$\rho(x,y):=\sum_{k=1}^\infty\frac{1}{2^k}(\sup_{t\in[\frac{1}{k},T]}\|x(t)-y(t)\|_{H^{s}}\wedge 1)$$
if the following conditions hold:

(a) For any $k\in \mathbb{N}$ there exists a constant $C(k)$ such that
$$\sup_{x\in B}\sup_{t\in [\frac{1}{k},T]}\|x(t)\|_{H^{s+1}}< C(k).$$

(b) $$\lim_{\delta\rightarrow0}\sup_{x\in B}\sum_{k=1}^\infty\frac{1}{2^k}(\sup_{|r-t|\leq \delta,\frac{1}{k}\leq r,t\leq T}\|x(t)-x(r)\|_{H^{s}}\wedge 1)=0$$

\proof Consider a sequence $\{x_n\}$ in $B$. By (a), (b) and the Arzela-Ascoli theorem for any $k\in \mathbb{N}$ we can find subsequence $\{x^k_{n_j}\}$ which is a subsequence of $\{x^{k-1}_{n_j}\}$ converging in  $C([\frac{1}{k},T];H^s)$. Then by a diagonal argument the result follows.$\hfill\Box$

\vskip.10in

\section{Local existence and uniqueness}
In this section we consider the following stochastic equation with multiplicative noise in $H$
$$d\theta+(f(\theta)+A^\alpha \theta)dt=G(\theta) dW, \eqno(3.1)$$
$$\theta(0)=\theta_0,$$
where $\alpha\in (0,1]$ and $W(t)$ is a cylindrical Wiener process in a separable Hilbert space $K$ defined on a filtered probability space $(\Omega,\mathcal{F},(\mathcal{F}_t)_{t\geq0},P)$. Here $G$ is a measurable mapping from $H$ to $L_2(K,H)$ (i.e. = all Hilbert-Schimit operator form $K$ to $H$) and we assume that there exists some $s_0\geq 1$ such that $f\in \Phi(H^{s_0}; H)\cap \Phi(H^{s_0+1+\alpha};H^{s_0})$ and  for every $\varphi\in \bigcap_{l=1}^\infty H^l$, $\langle f(\cdot),\varphi\rangle $ is continuous from $H^{s_0}$ to $\mathbb{R}$. Furthermore, $\theta_0$ is an $H^{s_0}$-valued $\mathcal{F}_0$-measurable function on $\Omega$.

Assume that $f$ satisfies the following conditions:
\vskip.10in
(b.1) (Coercivity condition) For every $s\in[s_0,s_0+1]$ there exists some locally bounded function $\rho_1$ on $\mathbb{R}$ and $\varepsilon_0\in(0,1)$ so that for every $v\in H^{2s}$ the following is satisfied:
$$-\langle f(v),\Lambda^{2s}v\rangle\leq \rho_1(|\Lambda^{s_0}v|)|\Lambda^{s_0}v|^2+\varepsilon_0|\Lambda^{s+\alpha}v|^2.\eqno(3.2)$$

\th{Remark} By interpolation and Young's inequality (3.2) is equivalent to the following:
$$-\langle f(v),\Lambda^{2s}v\rangle\leq \rho_1(|\Lambda^{s_0}v|)|\Lambda^{s}v|^2+\varepsilon_0|\Lambda^{s+\alpha}v|^2.$$

\vskip.10in

(b.2) (Growth condition) For $s=s_0,s_0-1$ there exists some locally bounded function $\rho_2$ on $\mathbb{R}$  such that for every $v\in H^{s_0+1+\alpha}$ $$|\Lambda^{s}f(v)|\leq \rho_2(|\Lambda^{s_0}v|)|\Lambda^{s+1+\alpha}v|.\eqno(3.3)$$

(b.3) (Local monotonicity condition) There exists some locally bounded function $\rho_3,\rho_4$ on $\mathbb{R}$  and $\tilde{\varepsilon}_0\in(0,1)$ such that for  $v_1,v_2\in H^{s_0}$
$$-\langle v_1-v_2,f(v_1)-f(v_2)\rangle\leq (\rho_3(|\Lambda^{s_0}v_1|)+\rho_4(|\Lambda^{s_0}v_2|))|v_1-v_2|^2+\tilde{\varepsilon}_0|\Lambda^\alpha (v_1-v_2)|^2.\eqno(3.4)$$
Also assume that $G$ satisfies the following conditions:

(G.1) For all $s\in [s_0,s_0+1]$, $G$ is an operator from $H^s$ to  $L_2(K,H^s)$ and there exists some locally bounded function $\rho_5$ on $\mathbb{R}$  such that for all $v\in H^s$ $$\|\Lambda ^sG(v)\|_{L_2(K,H)}\leq \rho_5(|\Lambda^{s_0}v|)(|\Lambda^sv|+1).$$

(G.2) There exists some locally bounded function $\rho_6,\rho_7$ on $\mathbb{R}$  such that for  $v_1,v_2\in H^{s_0}$
$$\|G(v_1)-G(v_2)\|_{L_2(K,H)}\leq (\rho_6(|\Lambda^{s_0}v_1|)+\rho_7(|\Lambda^{s_0}v_2|))|v_1-v_2|.$$

 \vskip.10in

First we make precise the notions of local, maximal and global solutions of (3.1).

\th{Definition 3.1} Fix a stochastic basis $(\Omega,\mathcal{F},P,\mathcal{F}_t,W)$.

(i) A local strong solution of (3.1) is a pair $(\theta,\tau)$, where $\tau$ is an $(\mathcal{F}_t)$-stopping time and $\theta=(\theta(t))_{t\geq0}$ is a predictable $H^{s_0}$-valued process such that $\theta(\cdot\wedge\tau)\in L^2(\Omega; L^2_{loc}([0,\infty);H^{s_0+\alpha}))$, $$\theta(\cdot\wedge\tau)\in C([0,\infty);H^{s_0})\quad P-a.s.,$$
and for every $t\geq0, \varphi\in \bigcap_{l=1}^\infty H^l$
$$\langle \theta(t\wedge \tau),\varphi\rangle+\int_0^{t\wedge \tau}\langle f(\theta)+A^\alpha \theta,\varphi\rangle dt=\langle\theta_0,\varphi\rangle+\langle\int_0^{t\wedge \tau}G(\theta) dW,\varphi\rangle \quad P-a.s..$$

(ii) We say that local pathwise uniqueness holds if given any pair $(\theta^1,\tau^1), (\theta^2,\tau^2)$ of local strong solutions of (3.1) with the same initial condition, the following holds:
$$P[\theta^1(t)=\theta^2(t);\forall t\in[0,\tau^1\wedge\tau^2]]=1.$$

(iii) A maximal strong solution of (3.1) is a pair $((\theta^R)_{R\in \mathbb{N}},(\tau_R)_{R\in \mathbb{N}})$ such that for each $R\in \mathbb{N}$ the pair $(\theta^R,\tau_R)$ is a local strong solution, $\tau_R$ is increasing such that $\zeta:=\lim_{R\rightarrow\infty}\tau_R>0$ $P$-a.s. and
$$\sup_{t\in [0,\tau_R]}|\Lambda^{s_0}\theta^R(t)|\geq R \qquad P-a.s.\textrm{ on the set } [\zeta<\infty].$$
 \vskip.10in
\th{Remark} If local pathwise uniqueness holds, then on $[t<\zeta]$ we can define $\theta(t):=\theta^n(t)$ on $[t\leq \tau_R]$. In this case we denote the maximal solution by $(\theta,(\tau_R)_{R\in \mathbb{N}},\zeta)$. We note that in this case $\zeta$ does not depend on the sequences $(\theta^R)_{R\in \mathbb{N}}, (\tau_R)_{R\in \mathbb{N}}$. Indeed, suppose there is another maximal strong solution $((\tilde{\theta}^R)_{R\in \mathbb{N}},(\tilde{\tau}_R)_{R\in \mathbb{N}})$ with $\tilde{\zeta}:=\lim_{R\rightarrow\infty}\tilde{\tau}_R$. Then up to a $P$-zero set
 $$\aligned \{\tilde{\zeta}>\zeta\}&\subset \cup_{n\in \mathbb{N}}\{\tilde{\tau}_n>\zeta\}\subset \cup_{n\in \mathbb{N}}\cap_{m\in\mathbb{N}}\{\tilde{\tau}_n>\tau_m,\zeta<\infty\}\\&\subset  \cup_{n\in \mathbb{N}}\cap_{m\in\mathbb{N}}\{\sup_{t\in[0,\tilde{\tau}_n]}|\Lambda^{s_0}\theta(t)|\geq m,\zeta<\infty\}\\&\subset  \cup_{n\in \mathbb{N}}\{\sup_{t\in[0,\tilde{\tau}_n]}|\Lambda^{s_0}\theta(t)|=\infty,\zeta<\infty\}.\endaligned$$
 Since $\theta\in C([0,\tilde{\tau}_n];H^{s_0})$, it follows that $P[\tilde{\zeta}>\zeta]=0$. Therefore, in this case a
  maximal strong solution $(\theta,(\tau_R)_{R\in \mathbb{N}},\zeta)$ is said to be global if $\zeta=\infty$ $P$- almost surely.
 \vskip.10in
 Now we want to show the  existence and uniqueness of local solutions of (3.1)

 \th{Theorem 3.2} Assume that $f$ satisfies (b.1-b.3) and $G$ satisfies (G.1)(G.2). Fix a stochastic basis $(\Omega,\mathcal{F},P,\mathcal{F}_t,W)$. Assume that $\theta_0$ is an $H^{s_0}$-valued, $\mathcal{F}_0$-measurable random variable with $E|\Lambda^{s_0}\theta_0|^{2}<\infty$. Then we have local pathwise uniqueness for (3.1) and there exists a maximal strong solution $(\theta,(\tau_R)_{R\in\mathbb{N}},\zeta)$ of (3.1).

 \proof [Step 1](Existence of local martingale solutions to (3.1))

 First  we construct the martingale solution to the following equation:
 $$d\theta+(\chi_R(|\Lambda^{s_0}\theta|)f(\theta)+A^\alpha \theta)dt=\chi_R(|\Lambda^{s_0}\theta|)G(\theta) dW,\eqno(3.5)$$
$$\theta(0)=\theta_0\in H^{s_0}.$$
Here we fix $R>0$ and $\chi_R:[0,\infty)\rightarrow[0,1]$ is a $C^\infty$ smooth function such that
 $$\chi_R(x)=\left\{\begin{array}{ll}1&\ \ \ \ \textrm{ for } x\leq R,\\0&\ \ \ \ \textrm{ for } x>2R.\end{array}\right.$$

 Consider the Galerkin approximation $\theta_n$ of (3.5) on $H_n$:   $$d\theta_n+(\chi_R(|\Lambda^{s_0}\theta_n|)P_nf(\theta_n)+A^\alpha \theta_n)dt=\chi_R(|\Lambda^{s_0}\theta_n|)P_nG(\theta_n) dW,\eqno(3.6)$$
$$\theta_n(0)=P_n\theta_0.$$Here $P_n$ is the projection operator onto $H_n$. Then by [PR07, Theorem 3.1.1] (b.2),(b.3), (G.1), (G.2) there exists a unique global solution $\theta_n$ to (3.6). By It\^{o}'s formula we obtain that for $s\geq s_0,$
$$  \aligned d|\Lambda^{s}\theta_n|^2+2|\Lambda^{s+\alpha}\theta_n|^2dt\leq& -2\chi_R(|\Lambda^{s_0}\theta_n|)\langle f(\theta_n),\Lambda^{2s}\theta_n\rangle dt+\chi_R(|\Lambda^{s_0}\theta_n|)\|\Lambda ^{s}P_nG(\theta_n)\|^2_{L_2(K,H)}dt\\&+2\langle \Lambda^{s}\theta_n, \chi_R(|\Lambda^{s_0}\theta_n|)\Lambda^{s}G(\theta_n)dW\rangle.\endaligned $$
By (b.1) and (G.1) we get for $s\in[s_0,s_0+1]$
$$  \aligned d|\Lambda^{s}\theta_n|^2+2|\Lambda^{s+\alpha}\theta_n|^2dt\leq & [C\chi_R(|\Lambda^{s_0}\theta_n|)|\Lambda^{s}\theta_n|^2+2\varepsilon_0|\Lambda^{s+\alpha}\theta_n|^2+C] dt\\&+2\chi_R(|\Lambda^{s_0}\theta_n|)\langle \Lambda^{s}\theta_n, \Lambda^{s}G(\theta_n)dW\rangle.\endaligned \eqno(3.7)$$
Then taking $s=s_0$ in (3.7) and by   the BDG-inequality for $p=1$, we have for any $T>0$
$$\aligned &E\sup_{t\in[0,T]}|\Lambda^{s_0}\theta_n(t)|^2+E\int_0^T|\Lambda^{s_0+\alpha}\theta_n(t)|^2dt\\\leq& CE|\Lambda^{s_0}\theta_0|^2+ CT+CE(\int_0^T|\Lambda^{s_0}\theta_n|^{2}\chi_R(|\Lambda^{s_0}\theta_n|)\|\Lambda^{s_0}G(\theta_n)\|^2_{L_2(K,H)}dt)^{1/2}.\endaligned.$$
Thus by (G.1) we get that for every $T>0$
$$E\sup_{t\in[0,T]}|\Lambda^{s_0}\theta_n(t)|^2+E\int_0^T|\Lambda^{s_0+\alpha}\theta_n(t)|^2dt\leq C_T,\eqno(3.8)$$
where $C_T$ is a constant independent of $n$.
Moreover, by (3.7), we obtain for $s\in[s_0,s_0+1]$, $0\leq t_0\leq t\leq T$,
$$E|\Lambda^{s}\theta_n(t)|^2+\int_{t_0}^tE|\Lambda^{s+\alpha}\theta_n(t)|^2dt\leq CE|\Lambda^{s}\theta_n(t_0)|^2+C\int_{t_0}^t(E(\chi_R(|\Lambda^{s_0}\theta_n(r)|)|\Lambda^{s}\theta_n(r)|^2)+1)dr.\eqno(3.9)$$
Now fix any $T$ and let $\varepsilon \in(0,T)$. By (3.8) we have that
$$\int_{\varepsilon/2}^\varepsilon E|\Lambda^{s_0+\alpha}\theta_n(t)|^2dt\leq C_T \qquad\textrm{ for all } n\in \mathbb{N}.\eqno(3.10)$$
Thus we can find $\tau_{\varepsilon,n}\in [\frac{\varepsilon}{2},\varepsilon]$ such that
$$\frac{\varepsilon}{2}E|\Lambda^{s_0+\alpha}\theta_n(\tau_{\varepsilon,n})|^2\leq C_T\qquad\textrm{ for all } n\in \mathbb{N}.\eqno(3.11)$$
Moreover, from (3.9) with $s=s_0+\alpha$ and $t=\varepsilon,t_0=\tau_{\varepsilon,n}$ we find that
$$E|\Lambda^{s_0+\alpha}\theta_n(\varepsilon)|^2\leq CE|\Lambda^{s_0+\alpha}\theta_n(\tau_{\varepsilon,n})|^2+C\int_{\varepsilon/2}^\varepsilon(E|\Lambda^{s_0+\alpha}\theta_n(r)|^2+1)dr\leq C(\frac{2}{\varepsilon}C_T+C_T+T),$$
where in the second inequality we used (3.10) and (3.11). Then taking $s=s_0+\alpha$ in (3.7) by the BDG-inequality, we have $$\aligned &E\sup_{t\in[\varepsilon,T]}|\Lambda^{s_0+\alpha}\theta_n(t)|^2+E\int_\varepsilon^T|\Lambda^{s_0+2\alpha}\theta_n(t)|^2dt\\\leq& CE|\Lambda^{s_0+\alpha}\theta_n(\varepsilon)|^2+\int_\varepsilon^T(CE|\Lambda^{s_0+\alpha}\theta_n|^2+C)dt\\&+CE(\int_\varepsilon^T|\Lambda^{s_0+\alpha}\theta_n|^{2}\chi_R(|\Lambda^{s_0}\theta_n|)\|\Lambda^{s_0+\alpha}G(\theta_n)\|^2_{L_2(K,H)}dt)^{1/2}\\\leq& CE|\Lambda^{s_0+\alpha}\theta_n(\varepsilon)|^2+\int_\varepsilon^T(CE\sup_{t\in[\varepsilon,r]}|\Lambda^{s_0+\alpha}\theta_n(t)|^2+C)dr+\frac{1}{2}E\sup_{t\in[\varepsilon,T]}|\Lambda^{s_0+\alpha}\theta_n(t)|^2.\endaligned.$$
Here in the second inequality we used (G.1) and Young's inequality. Then Gronwall's lemma yields that for some constant $C(\varepsilon)$
$$E\sup_{t\in[\varepsilon,T]}|\Lambda^{s_0+\alpha}\theta_n(t)|^2+E\int_\varepsilon^T|\Lambda^{s_0+2\alpha}\theta_n(t)|^2dt\leq C(\varepsilon)\qquad\textrm{ for all } n\in \mathbb{N}.$$
Then by a boot-strapping argument after finitely many steps we obtain that for some constant $C_1(\varepsilon)$
$$E\sup_{t\in[\varepsilon,T]}|\Lambda^{s_0+1}\theta_n(t)|^2+E\int_\varepsilon^T|\Lambda^{s_0+1+\alpha}\theta_n(t)|^2dt\leq C_1(\varepsilon)\qquad\textrm{ for all } n\in \mathbb{N}.\eqno(3.12)$$
 Now we prove that the family $\mathcal{L}(\theta_n),n\in \mathbb{N}$, is tight on $C([\frac{1}{k},T];H^{s_0})$, for all $k\in \mathbb{N}$. By (3.12) for each $t\in [\frac{1}{k},T]$, the laws of $\theta_n(t), n\in\mathbb{N},$ are tight on $H^{s_0}$. Then by Aldous' criterion in [Al78] it suffices to check that for all stopping times $\frac{1}{k}\leq\tau_n\leq T$ and $\eta_n\in \mathbb{R}$ with $\eta_n\rightarrow0,$ $$\lim_nE\|\theta_n(\tau_n+\eta_n)-\theta_n(\tau_n)\|_{H^{s_0}}=0.\eqno(3.13)$$
We have
$$\theta_n(\tau_n+\eta_n)-\theta_n(\tau_n)=-\int_{\tau_n}^{\tau_n+\eta_n}(\chi_R(|\Lambda^{s_0}\theta_n|)P_nf(\theta_n)+A^\alpha \theta_n)dt+\int_{\tau_n}^{\tau_n+\eta_n}\chi_R(|\Lambda^{s_0}\theta_n|)P_nG(\theta_n) dW.$$
By (3.12) we have that for large $n$
$$E\|\int_{\tau_n}^{\tau_n+\eta_n}A^\alpha\theta_n(t)dt\|_{H^{s_0}}\leq C\eta_n^{1/2}(E\int_{1/k}^{T+1}|\Lambda^{s_0+2\alpha}\theta_n(t)|^2dt)^{1/2}\rightarrow0, \textrm{ as } \eta_n\rightarrow0.$$
And by (b.2) and (3.12) we obtain that
$$\aligned &E\|\int_{\tau_n}^{\tau_n+\eta_n}(\chi_R(|\Lambda^{s_0}\theta_n|)P_nf(\theta_n)dt\|_{H^{s_0}}\\\leq& CE\int_{\tau_n}^{\tau_n+\eta_n}|\Lambda^{s_0+1+\alpha}\theta_n|dt\\\leq& C\eta_n^{1/2}(E\int_{1/k}^{T+1}|\Lambda^{s_0+1+\alpha}\theta_n|^2dt)^{1/2}\rightarrow0, \textrm{ as } \eta_n\rightarrow0.\endaligned$$
Also by (G.1) we have $$\aligned & E\|\int_{\tau_n}^{\tau_n+\eta_n}\chi_R(|\Lambda^{s_0}\theta_n|)P_nG(\theta_n)dW\|^2_{H^{s_0}}\leq C E\int_{\tau_n}^{\tau_n+\eta_n}\chi_R(|\Lambda^{s_0}\theta_n|)\|\Lambda^{s_0}G(\theta_n)\|_{L_2(K,H)}^2dt\\&\leq C\eta_n\rightarrow 0\textrm{ as } \eta_n\rightarrow0.\endaligned$$
Thus (3.13) follows, which implies the tightness of $\mathcal{L}(\theta_n),n\in\mathbb{N},$ in $C([\frac{1}{k},T];H^{s_0})$. This yields that for each $\eta>0$
$$\lim_{\delta\rightarrow0}\sup_nP(\sup_{|r-t|\leq \delta,\frac{1}{k}\leq r,t\leq T}\|\theta_n(t)-\theta_n(r)\|_{H^{s_0}}>\eta)=0,$$
which implies that for each $\eta>0$
$$\lim_{\delta\rightarrow0}\sup_nP\left[\sum_{k=1}^\infty\frac{1}{2^k}(\sup_{|r-t|\leq \delta,\frac{1}{k}\leq r,t\leq T}\|\theta_n(t)-\theta_n(r)\|_{H^{s_0}}\wedge 1)>\eta\right]=0.$$
Then for every $\varepsilon_1>0$ there exists a sequence $\delta_j,j\in \mathbb{N},$ such that $$\sup_nP\left[\sum_{k=1}^\infty\frac{1}{2^k}(\sup_{|r-t|\leq \delta_j,\frac{1}{k}\leq r,t\leq T}\|\theta_n(t)-\theta_n(r)\|_{H^{s_0}}\wedge 1)>\frac{1}{j}\right]<\frac{\varepsilon_1}{2^{j+1}}.$$
Let  $$\aligned B=&\big(\cap_k\{x\in C([1/k,T];H^{s_0}):\sup_{t\in[\frac{1}{k},T]}\|x(t)\|_{H^{s_0+1}}\leq C_1(\frac{1}{k})^{1/2}(k+N)\}\big)\\&\cap \big(\cap_j\{x\in C([1/k,T];H^{s_0}):\sum_{k=1}^\infty\frac{1}{2^k}(\sup_{|r-t|\leq \delta_j,\frac{1}{k}\leq s,t\leq T}\|x(t)-x(r)\|_{H^{s_0}}\wedge 1)\leq\frac{1}{j}\}\big).\endaligned$$
Here $C_1(\frac{1}{k})$ is the constant in (3.12) with $\varepsilon=\frac{1}{k}$.
Thus by Lemma 2.1 $B$ is a relatively compact set in $C((0,T];H^{s_0})$ and by (3.12) we can choose $N$ large enough such that $P(\theta_n\in B^c)<\varepsilon_1$ for all $n\in \mathbb{N}$, which implies that $\mathcal{L}(\theta_{n}),n\in\mathbb{N},$ are tight in $C((0,T];H^{s_0})$. By (3.8), a similar calculation as in the proof of (3.13) and Aldous' criterion in [Al78] we also obtain that $\mathcal{L}(\theta_{n}),n\in\mathbb{N}$, are tight in $C([0,T];H^{s_0-1})$, since $\theta_0\in L^2(\Omega;H^{s_0})$.

Therefore, we   find a subsequence, still denoted by $\theta_n$, such that $\mathcal{L}(\theta_n)$ converges weakly in $$C((0,T];H^{s_0})\cap C([0,T]; H^{s_0-1}).$$
By Skorohod's embedding theorem, there exist a stochastic basis $(\Omega^1,\mathcal{F}^1,\{\mathcal{F}^1_t\}_{t\geq0}, P^1)$ and, on this basis, $C((0,T];H^{s_0})\cap C([0,T];H^{s_0-1})$-valued random variables $\theta^1,\theta_n^1,n\geq1$, such that $\theta_n^1$ has the same law as $\theta_{n}$ on $C((0,T];H^{s_0})\cap C([0,T];H^{s_0-1})$, and $\theta_n^1\rightarrow \theta^1$ in $C((0,T];H^{s_0})\cap C([0,T];H^{s_0-1}),P^1$ -a.s. For $\theta_n^1$ we also have (3.8), thus for $\theta^1$, (3.8) is also satisfied by Fatou's lemma.

 For each $n\geq1$, define the $H_n$-valued process $$M_n^1(t):=\theta_n^1(t)-P_n\theta_0+\int_0^t(\chi_R(|\Lambda^{s_0}\theta^1_n|)P_nf(\theta^1_n)+A^\alpha \theta^1_n)dr.$$
In fact, $M_n^1$ is a square integrable martingale with respect to the filtration $$\{\mathcal{G}_n^1\}_t=\sigma\{\theta_n^1(r),r\leq t\}$$ with quadratic variation process $$\langle M_n^1\rangle_t=\int_0^t\chi_R(|\Lambda^{s_0}\theta_n^1(r)|)^2P_nG(\theta_n^1(r))G(\theta_n^1(r))^*P_ndr.$$ For all $r\leq t\in[0,T]$, all bounded continuous functions $\phi$ on  $C([0,r];H^{s_0-1})$, and all $v\in \cap_{l=1}^\infty H^l$, we have
$$E^1(\langle M_n^1(t)-M_n^1(r),v\rangle\phi(\theta_n^1|_{[0,r]}))=0$$ and $$\aligned &E^1((\langle M_n^1(t),v\rangle^2-\langle M_n^1(r),v\rangle^2-\int_r^t \chi_R(|\Lambda^{s_0}\theta_n^1|)^2\|G(\theta_n^1)^*P_nv\|^2_Kdr)\phi(\theta_n^1|_{[0,r]}))=0.\endaligned$$
By the BDG-inequality we have for any $p\geq 2$
$$\sup_n E^1|\langle M_n^1(t),v\rangle|^{2p}\leq C\sup_nE^1(\int_0^t\chi_R(|\Lambda^{s_0}\theta_n^1|)^2\|G(\theta_n^1)^*P_nv \|_K^2dr)^p<\infty.\eqno(3.14)$$
By (b.2), we have
$$\aligned &E^1\int_0^t|\chi_R(|\Lambda^{s_0}\theta_n^1|)\langle P_nf(\theta_n^1),v\rangle-\chi_R(|\Lambda^{s_0}\theta^1|)\langle f(\theta^1),v\rangle|dr\\\leq&E^1\int_\varepsilon^t|\chi_R(|\Lambda^{s_0}\theta_n^1|)\langle P_nf(\theta_n^1),v\rangle-\chi_R(|\Lambda^{s_0}\theta^1|)\langle f(\theta^1),v\rangle|dr\\&+\int^\varepsilon_0|\chi_R(|\Lambda^{s_0}\theta_n^1|)\langle P_nf(\theta_n^1),v\rangle-\chi_R(|\Lambda^{s_0}\theta^1|)\langle f(\theta^1),v\rangle|dr\\\leq&E^1\int_\varepsilon^t|\chi_R(|\Lambda^{s_0}\theta_n^1|)\langle f(\theta_n^1),v\rangle-\chi_R(|\Lambda^{s_0}\theta^1|)\langle f(\theta^1),v\rangle|dr\\&+E^1\int_\varepsilon^t|\chi_R(|\Lambda^{s_0}\theta_n^1|)\langle f(\theta_n^1),(P_n-I)v\rangle| dr\\&
+C|v|E^1\int^\varepsilon_0(\chi_R(|\Lambda^{s_0}\theta^1|)|\Lambda^{s_0-1}f(\theta^1)|+\chi_R(|\Lambda^{s_0}\theta_n^1|)|\Lambda^{s_0-1}f(\theta^1_n)|)dr\\\leq&E^1|\int_\varepsilon^t\chi_R(|\Lambda^{s_0}\theta_n^1|)\langle f(\theta_n^1),v\rangle-\chi_R(|\Lambda^{s_0}\theta^1|)\langle f(\theta^1),v\rangle|dr\\&+CE^1\int_\varepsilon^t|\Lambda^{s_0+\alpha}\theta_n^1||(P_n-I)v| dr\\&
+C|v|E^1\int^\varepsilon_0(|\Lambda^{s_0+\alpha}\theta^1|+|\Lambda^{s_0+\alpha}\theta_n^1|)dr.\endaligned$$
Letting first $n\rightarrow\infty$ and then $\varepsilon\rightarrow0$, by the continuity of $\langle f(\cdot),v\rangle$ on $H^{s_0}$ and (3.8) the latter converges to zero. This and (3.8) implies that
$$\lim_{n\rightarrow\infty}E^1|\langle M_n^1(t)-M^1(t),v\rangle|=0,$$and hence by (3.14)$$\lim_{n\rightarrow\infty}E^1|\langle M_n^1(t)-M^1(t),v\rangle|^2=0,$$where $$M^1(t):=\theta^1(t)-\theta^1(0)+\int_0^t(\chi_R(|\Lambda^{s_0}\theta^1|)f(\theta^1)+A^\alpha \theta^1)dr.$$
Taking the limit we obtain that for all $r\leq t\in[0,T]$, all bounded continuous functions $\phi$ on  $C([0,r];H^{s_0-1})$, and  $v\in \cap_{l=1}^\infty H^l$,
 $$E^1(\langle M^1(t)-M^1(r),v\rangle\phi(\theta^1|_{[0,r]}))=0.$$
and similarly by (G.1), (G.2) and (3.8) we obtain for $r\leq t\in[0,T]$,
$$E^1((\langle M^1(t),v\rangle^2-\langle M^1(r),v\rangle^2-\int_r^t\chi_R(|\Lambda^{s_0}\theta|)^2\| G(\theta)^*v\|_K^2 dr)\phi(\theta^1|_{[0,r]}))=0.$$
Thus by the martingale representation theorem (cf. [DZ92, Theorem 8.2], [O05, Theorem 2]) there exists a stochastic basis $(\tilde{\Omega},\tilde{P},\tilde{\mathcal{F}},(\tilde{\mathcal{F}}_t)_{t\geq0})$, a cylindrical Wiener process $\tilde{W}$ and an $(\tilde{\mathcal{F}}_t)$-adapted process $\tilde{\theta}$ with paths in $C([0,T]; H^{s_0-1})\cap C((0,T];H^{s_0})$ such that $\tilde{\theta}$ satisfies (3.5) with $W$ replaced by $\tilde{W}$ and $\tilde{\theta}(0)$ has the same distribution as $\theta_0$.

Now we want to show that $\tilde{\theta}(t)$ converges to $\tilde{\theta}(0)$ strongly in $H^{s_0}$ as $t\rightarrow0$. Since (3.8) is also satisfied by $\tilde{\theta}$, we have $\tilde{\theta}\in C([0,T];H^{s_0-1})\cap L^\infty([0,T];H^{s_0})$ which implies that $\tilde{\theta}$ is weakly continuous in $H^{s_0}$.
Consequently $$|\Lambda^{s_0}\tilde{\theta}(0)|\leq \liminf_{t\rightarrow0}|\Lambda^{s_0}\tilde{\theta}(t)|.\eqno(3.15)$$
Taking $s=s_0$ in (3.7) and integrating from $0$ to $t$ by the BDG-inequality with $p=1$ we have that
$$\aligned E\sup_{r\in[0,t]}|\Lambda^{s_0}\theta_n(r)|^2\leq E|\Lambda^{s_0}\theta_0|^2+C\int_{0}^t(E(\chi_R(|\Lambda^{s_0}\theta_n(r)|)|\Lambda^{s_0}\theta_n(r)|^2)+1)dr\\+CE(\int_0^t\chi_R(|\Lambda^{s_0}\theta_n(r)|)|\Lambda^{s_0}\theta_n(r)|^2(|\Lambda^{s_0}\theta_n(r)|^2+1)dr)^{1/2}.\endaligned$$
Thus, we get $$\tilde{E}\sup_{r\in[0,t]}|\Lambda^{s_0}\tilde{\theta}(r)|^2\leq \liminf_{n\rightarrow\infty}E\sup_{r\in[0,t]}|\Lambda^{s_0}\theta_n(r)|^2\leq E|\Lambda^{s_0}\theta_0|^2+Ct+Ct^{1/2},$$
which implies that
$$\tilde{E}\limsup_{t\rightarrow0}|\Lambda^{s_0}\tilde{\theta}(t)|^2\leq \tilde{E}|\Lambda^{s_0}\tilde{\theta}(0)|^2.\eqno(3.16)$$
Combining (3.15) and (3.16) we have that
$$\tilde{E}|\Lambda^{s_0}\tilde{\theta}(0)|^2\leq\tilde{E}\liminf_{t\rightarrow0}|\Lambda^{s_0}\tilde{\theta}(t)|^2\leq\tilde{E}\limsup_{t\rightarrow0}|\Lambda^{s_0}\tilde{\theta}(t)|^2\leq\tilde{E}|\Lambda^{s_0}\tilde{\theta}(0)|^2,$$
which implies by (3.15) that $$|\Lambda^{s_0}\tilde{\theta}(0)|^2=\lim_{t\rightarrow0}|\Lambda^{s_0}\tilde{\theta}(t)|^2.$$ This equality combined with the weak convergence implies that $\tilde{\theta}(t)$ converges to $\tilde{\theta}(0)$ strongly in $H^{s_0}$ as $t\rightarrow0$ which implies that  $\tilde{\theta}\in C([0,T];H^{s_0})$ for every $T>0$.

Define the stopping time
$$\tau_R:=\inf\{t\geq0:|\Lambda^{s_0}\tilde{\theta}(t)|\geq R\}.$$
Then $(\tilde{\theta},\tau_R)$ is a local martingale solution of (3.1) such that $\tilde{\theta}(\cdot\wedge \tau_R)\in C([0,\infty);H^{s_0})$ $P$-a.s and  $\tilde{\theta}(\cdot\wedge \tau_R)\in L^2(\Omega,L^2_{loc}([0,\infty),H^{s_0+\alpha}))$. Now define $\theta(t):=\tilde{\theta}(t\wedge \tau_R)$ . By the continuity of $\tilde{\theta}$ we have $$\tau_R=\inf\{t\geq0:|\Lambda^{s_0}\theta(t)|\geq R\}.\eqno(3.17)$$
Thus $\theta$ satisfies the following equation for any $\varphi\in \cap_{l=1}^\infty H^l$
 $$\langle \theta(t),\varphi\rangle+\int_0^{t\wedge\tau_R}\langle f(\theta)+A^\alpha \theta,\varphi\rangle dt=\langle\theta(0),\varphi\rangle+\langle\int_0^{t\wedge \tau_R} G(\theta) d\tilde{W},\varphi\rangle\quad \textrm{for all } t\geq0,\eqno(3.18)$$
 with $\tau_R$ given by (3.17).
Here $\tilde{W}$ is the cylindrical Wiener process from the martingale representation theorem.

[Step 2] (Pathwise uniqueness of (3.18)) Now we want to prove the pathwise uniqueness of solutions for (3.18) in $L^2(\Omega; L_{loc}^2([0,\infty);H^{s_0+\alpha}))$.  Assume that $\theta_1$ and $\theta_2$ are two solutions of (3.18) in $C([0,T];H^{s_0})$ with the same initial value $\theta_0$ on the same stochastic basis $(\Omega,\mathcal{F},P,\mathcal{F}_t,W)$. Take $v=\theta_1-\theta_2$ and we have for $t\in [0,\tau_R^1\wedge\tau_R^2]$
 $$v(t)+\int_0^t( f(\theta_1)- f(\theta_2)+A^\alpha v)dt=\int_0^t (G(\theta_1)-G(\theta_2)) dW,$$where $\tau_R^1$ and $\tau_R^2$ are given by (3.17) with $\theta$ replaced by $\theta_1$ and $\theta_2$ respectively. By It\^{o}'s formula, we obtain that for $t\in [0,\tau_R^1\wedge\tau_R^2]$
 $$\aligned d|v(t)|^2+2|\Lambda^\alpha v(t)|^2dt\leq&-2\langle v(t),f(\theta_1)- f(\theta_2)\rangle dt\\&+2\langle v(t),(G(\theta_1)-G(\theta_2)) dW\rangle+\|G(\theta_1)-G(\theta_2)\|^2_{L_2(K,H)}dt\endaligned$$
Then by (b.3)(G.2) and Young's inequality we deduce that for $t\in [0,\tau_R^1\wedge\tau_R^2]$
 $$\aligned d|v(t)|^2+2|\Lambda^\alpha v(t)|^2dt\leq&(2\tilde{\varepsilon}_0|\Lambda^\alpha v(t)|^2+C|v(t)|^2)dt\\&+2\langle v(t),(G(\theta_1)-G(\theta_2)) dW\rangle.\endaligned$$
 Thus, we have that for $t\in [0,\tau_R^1\wedge\tau_R^2]$
 $$ E|v(t)|^2 \leq C\int_0^tE|v(r)|^2dr.$$
Therefore, Gronwall's lemma yields that for $t\in [0,\tau_R^1\wedge\tau_R^2],$
$$E|v(t)|^2=0.$$
Hence we obtain that $P$-a.s $\theta_1(t)=\theta_2(t)$ for $t\in [0,\tau_R^1\wedge\tau_R^2]$. Then from (3.17) we get that $\tau_R^1=\tau_R^2$ $P$-a.s. Moreover by (3.18) we know that $\theta_i(t)=\theta_i(\tau_R^i)$ for $t> \tau_R^i, i=1,2$. Thus we obtain pathwise uniqueness of (3.18).

[Step 3](Existence of local strong solutions) Now fix a stochastic basis $(\Omega,\mathcal{F},P,\mathcal{F}_t,W)$.
By Steps 1, 2 and the Yamada-Watanabe Theorem in [Ku07, Theorem 3.14], we obtain that there exists a strong solution $\theta^R$ of (3.18) with $\tau_R=\inf\{t\geq0:|\Lambda^{s_0}\theta^R(t)|\geq R\}$ such that
 $\theta^R=\theta^R(\cdot\wedge\tau_R)\in L^2(\Omega; L^2_{loc}([0,\infty);H^{s_0+\alpha})$ and $\theta^R(\cdot\wedge\tau_R)\in C([0,\infty);H^{s_0})$ $P$-a.s.. This implies that $(\theta^R,\tau^R)$ is a local strong solution for (3.1) in the sense of Definition 3.1 (i).

 [Step 4](Local pathwise uniqueness)  Assume that $(\theta^1,\tau^1)$ and $(\theta^2,\tau^2)$ are two local strong solutions for (3.1). Then define the stopping time $$\tau_R=\inf\{t\geq0: |\Lambda^{s_0}\theta^1(t)|+|\Lambda^{s_0}\theta^2(t)|\geq R\}.$$
 By a modification of Step 2, we obtain $P[\theta^1(t)=\theta^2(t);\forall t\in[0,\tau^1\wedge\tau^2\wedge\tau_R]]=1$. Thus local pathwise uniqueness follows by taking the limit $R\rightarrow\infty$.

[Step 5](Existence of maximal strong solutions) Now take $R\in \mathbb{N}$ and define
$\zeta:=\lim_{R\rightarrow\infty}\tau_R$. Since $\theta_0\in L^2(\Omega;H^{s_0})$, we have $\zeta>0$ $P$-a.s.. Define $\theta(t):=\theta^R(t)$ on $[t\leq\tau_R]$. Hence obviously, $(\theta,(\tau_R)_{R\in\mathbb{N}},\zeta)$ is a maximal strong solution of (3.1).
$\hfill\Box$

\section{ Non-explosion of solutions driven by multiplicative noise }
In this section we consider the stochastic equation with linear multiplicative noise
$$d\theta+(f(\theta)+A^\alpha \theta)dt=\beta\theta dW, \eqno(4.1)$$
$$\theta(0)=\theta_0,$$
where in this case $\beta\in \mathbb{R}, \alpha\in (0,1]$ and $W$ is a single 1D Brownian motion.
It is easy to show that $G(\theta)f_i=\delta_{1i}\beta\theta$ satisfies (G.1), (G.2). Here $\{f_i\}$ is an orthonormal basis on $K$. Thus the result in Theorem 3.2 can be applied here. In this section we additionally assume the following condition which is stronger than (b.1):
 \vskip.10in

(b.1')  For every $s\in[s_0,s_0+1]$ there exists some strictly increasing continuous function $\rho_1$ on $\mathbb{R}$ satisfying $\lim_{R\rightarrow\infty}\rho_1(R)=\infty$ and $\varepsilon_0\in(0,1)$ so that for every $a\in(0,\infty)$ and every $v\in H^{2s}$
$$-\langle a^{-1}f(av),\Lambda^{2s}v\rangle\leq \rho_1(a|\Lambda^{s_0}v|)|\Lambda^{s_0}v|^2+\varepsilon_0|\Lambda^{s+\alpha}v|^2.$$
 \vskip.10in
Now we transform (4.1) to a random PDE. Consider the stochastic process
$$\gamma(t)=e^{-\beta W_t},t\geq0$$
Define $v:=\gamma \theta$ and by It\^{o}'s formula we have
$$\partial_t v+\frac{\beta^2}{2}v+\gamma f(\gamma^{-1}v)+A^\alpha v=0,\eqno(4.2)$$
$$v(0)=\theta_0.$$

\th{Theorem 4.1}  Assume that $f$ satisfies (b.1'), (b.2), and (b.3). Then there exists a map $\kappa:[1,\infty)\times (\beta^2_0,\infty)\rightarrow [0,\infty)$ defined by $\kappa(R,\beta^2):=\rho_1^{-1}(\frac{\beta^2}{4})^2\frac{1}{R}$, where $\rho_1$ is as in (b.1') and $\beta_0^2$ is such that $\rho^{-1}_1(\beta_0^2/4)=0$, satisfying
$$\lim_{\beta^2\rightarrow\infty}\kappa(R,\beta^2)=\infty,$$
 such that whenever for $\beta^2>\beta_0^2$
$$|\Lambda^{s_0}\theta_0|^2\leq\kappa(R,\beta^2), \qquad P-a.s.,$$
then $$P(\zeta=\infty)\geq1-\frac{1}{R^{1/4}}$$
and $$P(\lim_{t\rightarrow\infty}|\Lambda^{s_0}\theta(t)|=0)\geq1-\frac{1}{R^{1/8}}.$$
In particular, for every $\varepsilon>0$ and any given deterministic initial condition, the probability that solutions corresponding to sufficiently large $|\beta|$ never blow up, is greater than $1-\varepsilon$.

\proof By Theorem 3.2 we obtain a maximal strong solution $(\theta,(\tau_R)_{R\in\mathbb{N}},\zeta)$ for (4.1). Define $v:=\gamma\theta$. In order to get a good estimate for $v$, we first prove that the solutions of the Galerkin approximations converge to $v$. Fix $\omega\in \Omega$ and consider the following Galerkin approximation  to the cutoff equation of (4.2):
$$\partial_t v_n+\frac{\beta^2}{2}v_n+\chi_R(|\Lambda^{s_0}(\gamma^{-1}v_n)|)P_n \gamma f(\gamma^{-1}v_n)+A^\alpha v_n=0,\eqno(4.3)$$
$$v(0)=P_n\theta_0,$$
where $P_n$ is the projection operator from (3.6).
Multiplying both sides of (4.3) by $A^{s}v_n$ and taking inner product in $L^2$ we obtain that for $s\in [s_0,s_0+1]$
$$\frac{d|\Lambda^{s}v_n|^2}{dt}+\beta^2|\Lambda^{s}v_n|^2+2|\Lambda^{s+\alpha}v_n|^2=-2\chi_R(|\Lambda^{s_0}(\gamma^{-1}v_n)|)\langle \gamma f(\gamma^{-1}v_n),\Lambda^{2s}v_n\rangle.$$
By (b.1') we get
$$\frac{d|\Lambda^{s}v_n|^2}{dt}+\beta^2|\Lambda^{s}v_n|^2+2|\Lambda^{s+\alpha}v_n|^2\leq 2\varepsilon_0|\Lambda^{s+\alpha}v_n|^2+2\chi_R(|\Lambda^{s_0}(\gamma^{-1}v_n)|)\rho_1(\gamma^{-1}|\Lambda^{s_0}v_n|)|\Lambda^sv_n|^2.\eqno(4.4)$$
Then taking $s=s_0$ in (4.4) by Gronwall's lemma we have that
$$\sup_{t\in[0,T]}|\Lambda^{s_0}v_n(t)|^2+\int_0^T|\Lambda^{{s_0}+\alpha}v_n|^2dt\leq C(\omega),$$
where $C(\omega)$ is a constant depending on $\omega$, but is independent of $n$.
By a similar argument as in the proof of Theorem 3.2, [Step 1] we obtain that$$\sup_{t\in[\varepsilon,T]}|\Lambda^{s_0+1}v_n(t)|^2+\int_\varepsilon^T|\Lambda^{s_0+1+\alpha}v_n|^2dt\leq C_2(\varepsilon,\omega).$$
Here $C_2(\varepsilon,\omega)$ is a constant depending on $\varepsilon$ and $\omega$. By a similar calculation as in the proof of Theorem 3.2, [Step 1] and the Arzela-Ascoli theorem, $\{v_n\}$ is compact in $C([0,T];H^{s_0-1})\cap C((0,T];H^{s_0})$ which implies that there exists a subsequence still denoted by $v_n$ converging to some $\tilde{v}$ in $C([0,T];H^{s_0-1})\cap C((0,T];H^{s_0})$. Also by (4.4) it is easy to obtain that $\tilde{v}\in C([0,T];H^{s_0})$. Now define $$\tilde{\tau}_R:=\inf\{t\geq0:|\Lambda^{s_0}(\gamma^{-1}\tilde{v}(t))|\geq R\}.$$
Then $\tilde{v}$ satisfies the following equation for any $\varphi\in \cap_{l=1}^\infty H^l$
$$\langle\tilde{v}(t\wedge\tilde{\tau}_R),\varphi\rangle-\langle\theta_0,\varphi\rangle+\frac{\beta^2}{2}\int_0^{t\wedge\tilde{\tau}_R}\langle\tilde{v},\varphi\rangle dr+\int_0^{t\wedge\tau_R}\langle \gamma f(\gamma^{-1}\tilde{v})+A^\alpha \tilde{v},\varphi\rangle dr=0\quad \textrm{ for all }t\geq0.\eqno(4.5)$$

 The $\omega$-wise uniqueness of (4.5) can be proved similarly to Theorem 3.2, [Step 2]. By It\^{o}'s formula,
 $v(t)$ is a solution to (4.5) which implies that $v(t)=\tilde{v}(t), t\in [0,\tau_R\wedge \tilde{\tau}_R], \tau_R=\tilde{\tau}_R$.

 Taking the limit in (4.4) with $s=s_0$ we obtain that for $0\leq r\leq t\leq\tau_R$ with $R$ large enough
$$|\Lambda^{s_0}v(t)|^2+\beta^2\int_r^t|\Lambda^{s_0}v|^2dt_1\leq |\Lambda^{s_0}v(r)|^2+2\int_r^t\rho_1(\gamma^{-1}|\Lambda^{s_0}v|)|\Lambda^{s_0}v|^2dt_1.$$
For the case that $\rho_1^{-1}(\frac{\beta^2}{4})>0$ we now define the stopping time
$$\sigma:=\inf\{t\geq0:\rho_1(\gamma^{-1}|\Lambda^{s_0}v|)> \frac{\beta^2}{4}\}=\inf\{t\geq0:|\Lambda^{s_0}\theta|> \rho_1^{-1}(\frac{\beta^2}{4})\}.$$
For $R$ large enough we have $\sigma\leq\tau_R$.
Then we obtain on $[0,\sigma]$,
$$\frac{d|\Lambda^{s_0}v|^2}{dt}+\frac{\beta^2}{2}|\Lambda^{s_0}v|^2\leq0.$$
Thus we have on $[0,\sigma]$, $$|\Lambda^{s_0}v(t)|^2\leq |\Lambda^{s_0}\theta_0|^2\exp(-\frac{\beta^2t}{2}).$$
Recalling that for our maximal solution $(\theta,(\tau_R)_{R\in\mathbb{N}},\zeta)$ we have $\theta=\gamma^{-1}v$, we deduce that
$$|\Lambda^{s_0}\theta(t)|^2\leq \gamma^{-2}|\Lambda^{s_0}\theta_0|^2\exp(-\frac{\beta^2t}{2}).\eqno(4.6)$$
Set $\rho(t):=\gamma^{-2}\exp(-\frac{\beta^2t}{2})=\exp(2\beta W_t-\frac{\beta^2t}{2}).$ Fix any $R\geq1$ and define the stopping time
$$\sigma_R:=\inf\{t\geq0:\rho(t)\geq R\}.$$
Then for  $t\in[0,\sigma\wedge \sigma_R)$ we have
$$|\Lambda^{s_0}\theta(t)|^2\leq R|\Lambda^{s_0}\theta_0|^2.\eqno(4.7)$$
Recall that $\kappa(R,\beta^2):=\rho_1^{-1}(\frac{\beta^2}{4})^2\frac{1}{R}$. So, $$\lim_{\beta^2\rightarrow\infty}\kappa(R,\beta^2)=\infty,$$
and if the initial datum satisfies $$|\Lambda^{s_0}\theta_0|^2\leq \kappa(R,\beta^2),$$ then by (4.7) for all $t\in [0,\sigma\wedge \sigma_R]$,
$$|\Lambda^{s_0}\theta(t)|^2\leq \rho_1^{-1}(\frac{\beta^2}{4})^2,$$ which implies $$|\Lambda^{s_0}\theta(t)|\leq \rho_1^{-1}(\frac{\beta^2}{4}).$$
Hence due to the definition of $\sigma$, $\sigma\wedge\sigma_R=\sigma_R$  and hence $\zeta\geq \tau_R\geq\sigma\geq\sigma\wedge\sigma_R=\sigma_R$. Therefore, the maximal strong solution $(\theta,(\tau_R)_{R\in\mathbb{N}},\zeta)$ is global in time on the set $\{\sigma_R=\infty\}$. By the martingale maximal inequality it follows ( see e.g. [GV12, Lemma 9.1]) that $$P(\sigma_R=\infty)\geq1-\frac{1}{R^{1/4}}.$$Thus the first assertion follows.
Now we prove the second result.
Set $\rho_0(t):=\gamma^{-2}\exp(-\frac{\beta^2t}{4})=\exp(2\beta W_t-\frac{\beta^2t}{4}).$ Fix any $R\geq1$ and define the stopping time
$$\sigma^0_R:=\inf\{t\geq0:\rho_0(t)\geq R\}.$$
Then by (4.6) on $[0,\sigma\wedge \sigma^0_R]$ we have
$$|\Lambda^{s_0}\theta(t)|^2\leq R|\Lambda^{s_0}\theta_0|^2\exp(-\frac{\beta^2t}{4}).$$
Thus, if the initial datum satisfies $$|\Lambda^{s_0}\theta_0|^2\leq \kappa(R,\beta^2),$$ then for all $t\in [0,\sigma\wedge \sigma^0_R]$,
$$|\Lambda^{s_0}\theta(t)|^2\leq \rho_1^{-1}(\frac{\beta^2}{4})^2\exp(-\frac{\beta^2t}{4}).\eqno(4.7)$$
Hence due to the definition of $\sigma$, $\sigma\wedge\sigma^0_R=\sigma^0_R.$ Hence $\zeta\geq\sigma\geq \sigma^0_R$. Therefore the maximal strong solution $(\theta,(\tau_R)_{R\in\mathbb{N}},\zeta)$ is global in time and $|\Lambda^{s_0}\theta(t)|\rightarrow0$ as $t\rightarrow\infty$ on the set $\{\sigma^0_R=\infty\}$ .  But again by the martingale maximal inequality we have $$P(\sigma^0_R=\infty)\geq1-\frac{1}{R^{1/8}}.$$
Thus the second assertion follows.
$\hfill\Box$

\section{Application to some examples}
In this section we describe some examples for (3.1) satisfying conditions (b.1'),(b.2),(b.3) imposed above. First, we recall some useful estimates which will be used later.

Let $O$ be a bounded open domain in $\mathbb{R}^d$ with smooth boundary and let $C^\infty_c(O)$  denote the set of all smooth functions from
$O$ to $\mathbb{R}$ with compact supports. For $p > 1$, let $L
^p(O)$ be the  $L^p$
-space in which the norm is denoted by $\|\cdot\|_{L^p}$. If $A$ is $-\Delta$ on the domain $O$ with Dirichlet or periodic boundary condition, then we have the following estimates which will be used later. For $s \geq0, p \in[1,\infty]$ we use $H_0^{s,p}(O)$ (or $H^{s,p}(O)$) to denote the Sobolev space os all $f\in H$ for which $\|\Lambda^sf\|_{L^p}$
is finite. Let us start with the following important product estimates from [Re95, Lemma A.4]:
\vskip.10in
\th{Lemma 5.1} Suppose that $s>0$ and $p\in (1,\infty)$. If $f,g\in C_c^\infty(O)$ or $C^\infty(\mathbb{T}^d)$, then
$$\|\Lambda^s(fg)\|_{L^p}\leq C(\|f\|_{L^{p_1}}\|\Lambda^sg\|_{L^{p_2}}+\|g\|_{L^{p_3}}\|\Lambda^sf\|_{L^{p_4}}),$$
with $p_i\in (1,\infty], i=1,...,4$ such that
$$\frac{1}{p}=\frac{1}{p_1}+\frac{1}{p_2}=\frac{1}{p_3}+\frac{1}{p_4}.$$
\vskip.10in

We shall also use the following standard Sobolev inequality (cf. [St70, Chapter V]):
\vskip.10in
\th{Lemma 5.2} Suppose that $q>1, p\in [q,\infty)$ and
$$\frac{1}{p}+\frac{\sigma}{d}=\frac{1}{q}.$$
Suppose that $\Lambda^\sigma f\in L^q$, then $f\in L^p$ and there is a constant $C\geq 0$ such that
$$\|f\|_{L^p}\leq C\|\Lambda^\sigma f\|_{L^q}.$$
\vskip.10in
The following commutator estimates from [Ju04, Lemma 3.1] are very important for later use.

\th{Lemma 5.3 (Commutator Estimates)} Suppose that $s>0$ and $p\in (1,\infty)$. If $f,g\in C^\infty_c(O)$ or $C^\infty(\mathbb{T}^d)$, then
$$\|\Lambda^s(fg)-f\Lambda^s(g)\|_{L^p}\leq C(\|\nabla f\|_{L^{p_1}}\|\Lambda^{s-1}g\|_{L^{p_2}}+\|g\|_{L^{p_3}}\|\Lambda^sf\|_{L^{p_4}}),$$
with $p_i\in (1,\infty), i=1,...,4$ such that
$$\frac{1}{p}=\frac{1}{p_1}+\frac{1}{p_2}=\frac{1}{p_3}+\frac{1}{p_4}.$$
\vskip.10in
\th{Remark 5.4} Now we  give examples for $G$ satisfying (G.1) and (G.2) when $H=L^2(O)$ and $A$ is as in Section 2. Let $f_n, n\in \mathbb{N},$ be an ONB of $K$. For $s\in\mathbb{R}^+$, we write $s=[s]+\{s\}$ with $[s]\in \mathbb{Z}, \{s\}\in (0,1)$. For $y\in K$
$$G(\theta)y=\sum_{k=1}^\infty b_kg(\theta)\langle y,f_k\rangle_K,\theta\in H^\alpha,$$
where $b_k\in C^\infty(\mathbb{T}^d)$ and $g:\mathbb{R}\mapsto \mathbb{R}$ is $C^\infty$ smooth.
Assume $s_0\geq \frac{d}{2}+1$,
$$(\sum_k |\Lambda^{s_0+1} b_k|^2)^{1/2}\leq M,$$
and that there exists $m_5>0$ such that for $k=1,...,[\frac{d}{2}]+[s_0]+2$,  $$|g^{(k)}(\xi)|\leq C(1+|\xi|^{m}),$$
where $g^{(k)}$ denotes the $k$-th derivative of $g$.

Now for $s\in[s_0,s_0+1], $ we have
$$\|\Lambda ^sG(\theta)\|_{L_2(K,H)}\leq M|\Lambda^sg(\theta)|\leq M|\Lambda^{\{s\}}D^{[s]}g(\theta)|.$$
 Since $$D^{[s]}g(\theta)= \sum_{\beta_1+...+\beta_\mu=[s]}C_\beta\theta^{(\beta_1)}\cdot...\cdot\theta^{(\beta_\mu)}g^{(\mu)}(\theta),$$
 where $\theta^{(\beta_i)}$ denotes the $\beta_i$-th derivative of $\theta$,
 we have$$\|\Lambda ^sG(\theta)\|_{L_2(K,H)}\leq M\sum_{\beta_1+...+\beta_\mu=[s]}C_\beta|\Lambda^{\{s\}}(\theta^{(\beta_1)}\cdot...\cdot\theta^{(\beta_\mu)}g^{(\mu)}(\theta))|.$$
 Now we estimate each term of the right hand side of the above inequality: for $\mu=1$ by Lemmas 5.1, 5.2 we have
 $$\aligned|\Lambda^{\{s\}}(\theta^{([s])}g^{(1)}(\theta))|\leq & C[|\Lambda^{\{s\}}\theta^{([s])}|\|g^{(1)}(\theta)\|_{L^\infty}+|\Lambda^{\{s\}}\theta^{([s])}|\|\Lambda^{\{s\}}g^{(1)}(\theta)\|_{L^p}]
 \\\leq & C[|\Lambda^{\{s\}}\theta^{([s])}|\|g^{(1)}(\theta)\|_{L^\infty}+|\Lambda^{\{s\}}\theta^{([s])}|\|\nabla\theta(|\theta|^m+1)\|_{L^p}]\\\leq&C[|\Lambda^{s_0}\theta|^{m+1}+1]|\Lambda^s\theta|.\endaligned$$
 Here $\frac{1}{p}=(\{s\}/d)\wedge(\frac{1}{2}-\varepsilon)$ with some $\varepsilon>0$ and we used $H^{s_0}\subset H_0^{1,p}(O)\subset H_0^{\{s\},p}(O)$ (or $H^{s_0}\subset H^{1,p}(O)\subset H^{\{s\},p}(O)$) and $H^{s_0}\subset L^\infty$ in the last inequality.
 For $\mu=2$, by Lemmas 5.1, 5.2 we have
 $$\aligned\sum_{\beta_1+\beta_2=[s]}|\Lambda^{\{s\}}(\theta^{(\beta_1)}\theta^{(\beta_2)}g^{(2)}(\theta))|\leq & C\sum_{\beta_1+\beta_2=[s]}|\Lambda^{\{s\}}(\theta^{(\beta_1)}\theta^{(\beta_2)})|\cdot[\|g^{(2)}(\theta)\|_{L^\infty}+
 \|\Lambda^{\{s\}}g^{(2)}(\theta)\|_{L^p}]\\\leq& C\sum_{\beta_1+\beta_2=[s]}|\Lambda^{\{s\}+\sigma_1}\theta^{(\beta_1)}||\Lambda^{\sigma_2}\theta^{(\beta_2)}|\cdot[\|g^{(2)}(\theta)\|_{L^\infty}+
 \|\nabla\theta(|\theta|^m+1)\|_{L^p}]\\\leq& C|\Lambda^{s}\theta||\Lambda^{\frac{d}{2}}\theta|\cdot[|\Lambda^{s_0}\theta|^{m+1}+1],\endaligned$$
 where $\frac{1}{p}=(\{s\}/d)\wedge(\frac{1}{2}-\varepsilon)$ with some $\varepsilon>0$ and $\sigma_1,\sigma_2>0,\beta_1+\sigma_1\leq[s], \beta_2+\sigma_2\leq s_0, \sigma_1+\sigma_2=\frac{d}{2}$. Here we used $H^{s_0}\subset H_0^{1,p}(O)\subset H_0^{\{s\},p}(O)$ (or $H^{s_0}\subset H^{1,p}(O)\subset H^{\{s\},p}(O)$) and $H^{s_0}\subset L^\infty$ in the last inequality.
The other terms can be estimated similarly and (G.1) follows. Furthermore,
$$\aligned \|G(\theta_1)-G(\theta_2)\|_{L_2(K,H)}\leq& M |g(\theta_1)-g(\theta_2)|\\\leq& M\int_0^1\|g^{(1)}(r\theta_1+(1-r)\theta_2)\|_{L^\infty}|\theta_1-\theta_2|dr\\\leq& C(|\Lambda^{s_0}\theta_1|^{m}+|\Lambda^{s_0}\theta_2|^{m}+1)|\theta_1-\theta_2|,\endaligned$$
hence (G.2) holds.

\vskip.10in
In the following subsections we consider the same situation as described in Section 3 and 4, but give concrete examples for $H, A$ and $f$ respectively. We fix a stochastic basis $(\Omega,\mathcal{F},(\mathcal{F}_t)_{t\geq0},P)$ and a cylindrical Wiener process $W$ on $H$.

\subsection { Stochastic fractional Burgers equation in the supercritical case}

 We consider the following stochastic fractional Burgers equation in $\mathbb{T}^1$:
$$d\theta+(\theta\nabla\theta+(-\Delta)^\alpha \theta)dt= G(\theta) dW. \eqno(5.1)$$

  The stochastic Burgers equation with $\alpha=1$ has received an extensive amount of
attention (see e.g. [DDT94], [DZ96] and the references therein). In these papers, the authors obtained existence and uniqueness of global solutions for $\alpha=1$. Recently to study the relation between nonlinear and dissipative phenomena, many researchers studied the case when $\alpha<1$ in the deterministic case.  A. Kiselev, F. Nazarov and R. Shterenberg obtained in [KNS08] that if $\alpha<\frac{1}{2}$, the solution to (5.1) in the deterministic case may  blow up in finite time. The  stochastic fractional Burgers equation driven by space-time white noise
for general parameter $\alpha\in(\frac{3}{4}, 1)$  has been studied in [BGD11], in which the authors obtain
the existence and uniqueness of mild solutions. Now we apply our Theorems 3.2 and 4.1 to (5.1).

Let
$$H:=\{f\in L^2(\mathbb{T}^1;\mathbb{R}):\int_{\mathbb{T}^1}f d\xi=0\}.$$
 $$Av:=-\Delta v\quad v\in D(A):= H^{2,2}(\mathbb{T}^1).$$
 Then $A$ is an unbounded positive definite  self-adjoint operator and   $A^{-1}$ is compact on $H$. Then $H^{s}, s\in\mathbb{R}$, are the classical Sobolev spaces on $\mathbb{T}^1$. We set
 $$f(v):=v\nabla v \quad v\in H^1.$$
We can easily find an $s_0\geq 1$ such that $f\in \Phi(H^{s_0};H)\cap \Phi(H^{s_0+1+\alpha};H^{s_0})$ and for every $\varphi\in \bigcap_{l=1}^\infty H^l$, $\langle f(\cdot),\varphi\rangle $ is continuous from $H^{s_0}$ to $\mathbb{R}$.
 Moreover, by [KNS08, Lemma 2.1], we obtain that for any positive constant $a$, $s>\frac{3}{2}-2\alpha$ and $v\in H^{s+\alpha}\cap H^{2s}$
$$|\langle a\partial_\xi(v^2),\Lambda^{2s}v\rangle|\leq Ca|\Lambda^{s+\alpha-\varepsilon_2}v|^2|\Lambda^qv|\leq \varepsilon_1|\Lambda^{s+\alpha}v|^2+C_2a^{1/\delta}|\Lambda^qv|^{2+\frac{1}{\delta}},$$
where $s+\alpha-\varepsilon_2>q>\frac{3}{2}-2(\alpha-\varepsilon_2)$ and $\delta:=\delta(s)=\frac{\varepsilon_2}{s+\alpha-q}$. Here we used the interpolation inequality  and Young's inequality in the last inequality.
Thus, (b.1') is satisfied if $s_0> \frac{3}{2}-2\alpha$ and $s_0\geq1$. By Lemmas 5.1 and 5.2 we have for $s=s_0,s_0-1$, $v\in H^{s_0+1+\alpha}$ that
$$|\Lambda^{s}(v\nabla v)|\leq C|\Lambda^{s+1+\alpha}v||\Lambda v|.$$
Then (b.2) is satisfied for $s_0\geq1$. Now we only need to prove (b.3).
By a similar calculation as in the proof of [KNS08, Theorem 2.8] we have that for $\theta_1,\theta_2\in H^{s_0}\cap H^{\frac{3}{2}-\alpha},$
$$\aligned |\langle \partial_\xi(\theta_1^2)-\partial_\xi(\theta_2^2),\theta_1-\theta_2\rangle|=&|2\langle \partial_\xi(\theta_1),(\theta_1-\theta_2)^2\rangle+2\langle\partial_\xi(\theta_1-\theta_2),\theta_2(\theta_1-\theta_2)\rangle|\\=&|2\langle \partial_\xi(\theta_1),(\theta_1-\theta_2)^2\rangle-2\langle\partial_\xi(\theta_1-\theta_2),(\theta_1-\theta_2)^2\rangle\\&+2\langle\partial_\xi(\theta_1-\theta_2),\theta_1(\theta_1-\theta_2)\rangle|
\\=&|\langle \partial_\xi(\theta_1),(\theta_1-\theta_2)^2\rangle|\\\leq& C|\Lambda^{(\frac{3}{2}-\alpha)\vee1}\theta_1||\theta_1-\theta_2||\Lambda^\alpha(\theta_1-\theta_2)|.\endaligned$$Thus (b.3) is satisfied if $s_0\geq (\frac{3}{2}-\alpha)\vee1$.
Here in the third equality we used the integration by parts and in the last inequality we used Lemma 5.2.

 Now Theorems 3.2 and 4.1 apply to give the following results:
\vskip.10in
\vskip.10in
\th{Theorem 5.5} Fix $0<\alpha<1$. Assume that  $G$ satisfies (G.1) and (G.2) with $s_0\geq (\frac{3}{2}-\alpha)\vee1$. Assume that $\theta_0$ is an $H^{s_0}$-valued, $\mathcal{F}_0$-measurable random variable with $E|\Lambda^{s_0}\theta_0|^2<\infty$.

(i) Then local pathwise uniqueness holds and there exists a maximal strong solution $(\theta,(\tau_R)_{R\in\mathbb{N}},\zeta)$ of (5.1).

(ii) Moreover, for $G(\theta)f_i=\delta_{1i}\beta\theta$  there exists a positive deterministic funtion $\kappa:[1,\infty)\times (0,\infty)\rightarrow (0,\infty)$ defined by $\kappa(R,\beta^2):=(\frac{\beta^2}{4C_1})^{2\delta(s_0)}\frac{1}{R}$, where $C_1,\delta(s_0)$ are as above, satisfying
$$\lim_{\beta^2\rightarrow\infty}\kappa(R,\beta^2)=\infty,$$
 such that whenever
$$|\Lambda^{s_0}\theta_0|^2\leq\kappa(R,\beta^2) \qquad P-a.s.,$$
then $$P(\zeta=\infty)\geq1-\frac{1}{R^{1/4}}$$
and $$P(\lim_{t\rightarrow\infty}|\Lambda^{s_0}\theta(t)|=0)\geq1-\frac{1}{R^{1/8}}.$$

\subsection{ Stochastic quasi-geostrophic equation in the supercritical case}

 We consider the following stochastic quasi-geostrophic equation in $\mathbb{T}^2$: for $0<\alpha<1$
$$d\theta+(u\cdot\nabla\theta+(-\Delta)^\alpha \theta)dt=G(\theta) dW, \eqno(5.2)$$
where $$u=(u_1,u_2)=(-R_2\theta,R_1\theta)=R^\bot\theta.$$

This equation is an important model in geophysical fluid dynamics. The case $\alpha = 1/2$
exhibits similar features (singularities) as the 3D Navier-Stokes equations and can therefore
serve as a model case for the latter. In the deterministic case, the global existence of weak solutions has been
obtained in [Re95] and one most remarkable result in [CV06] gives the existence of a classical
solution for $\alpha = 1/2$. In [KNV07] another very important result is proved, namely that solutions
for $\alpha = 1/2$ with periodic $C^\infty$ data remain $C^\infty$ for all times. The blow up or global regularity for $\alpha < 1/2$ remains an open problem for the  quasi-geostrophic equation. For more details we refer to [CCCGW12] and the reference therein.

The 2D stochastic quasi-geostrophic equation on $\mathbb{T}^2$
for general parameter $\alpha\in(0, 1)$  has been studied in [RZZ12], in which the authors obtain
the existence of martingale solutions for (5.2) for general parameter $\alpha\in(0, 1)$ and for both additive as well as multiplicative noise.

By the singular integral theory of Calder\'{o}n and Zygmund (cf [St70, Chapter 3]), for any $p\in(1,\infty)$, there is a constant $C=C(p)$, such that
$$\|R_j\theta\|_{L^p}\leq C(p)\|\theta\|_{L^p}.$$
Here $R_j$ is the $j$-th periodic Riesz transform.

 Now we apply Theorems 3.2 and 4.1 to (5.2).
Let $H:=\{f\in L^2(\mathbb{T}^2;\mathbb{R}),\int_{\mathbb{T}^2} fd\xi=0\}$ and
$$Av:=-\Delta v\textrm{,}\quad v\in D(A):= H^{2,2}(\mathbb{T}^2).$$
 Then $A$ is an unbounded positive definite  self-adjoint operator and   $A^{-1}$ is compact on $H$. Then $H^{s}, s\in\mathbb{R},$ are the classical Sobolev spaces on $\mathbb{T}^2$. Set
 $$f(v):=(-R_2v,R_1v)\cdot\nabla v=R^{\perp}v\cdot\nabla v\textrm{,} \quad v\in H^{1+\delta_0},$$ for some $\delta_0>0$
 with $R_jv$ being the $j$-th periodic Riesz transform.

We can easily find an $s_0>1$ such that $f\in \Phi(H^{s_0};H)\cap \Phi(H^{s_0+2+\alpha};H^{s_0+1})$ and for every $\varphi\in \bigcap_{l=1}^\infty H^l$, $\langle f(\cdot),\varphi\rangle $ is continuous from $H^{s_0}$ to $\mathbb{R}$.

Because $\langle u_v\cdot\nabla \Lambda^sv,\Lambda^{s}v\rangle=0$ for $v\in H^{2s}$ and by Lemma 5.3 , we obtain that for any positive constant $a$ and $s\geq 2-\alpha$, $v\in H^{2s}$,
$$\aligned |\langle au_v\cdot \nabla v,\Lambda^{2s}v\rangle|&=a|\langle \Lambda^{s}(u_v\cdot \nabla v)-u_v\cdot\nabla \Lambda^sv,\Lambda^{s}v\rangle|\\&\leq Ca|\Lambda^sv|\|\Lambda^sv\|_{L^q}\|\Lambda v\|_{L^p}\\&\leq Ca|\Lambda^sv||\Lambda^{s+\alpha}v||\Lambda^{2-\alpha} v|\\&\leq \varepsilon|\Lambda^{s+\alpha}v|^2+C_1a^{2/\delta}|\Lambda^{2-\alpha}v|^{2+\frac{2}{\delta}}.\endaligned$$
Here $u_v=R^\bot v, \frac{1}{p}+\frac{1}{q}=\frac{1}{2},\frac{1}{p}=\frac{\alpha}{2}, \delta:=\delta(s)=\frac{\alpha}{s+2\alpha-2}$ and we used $H^{\alpha}\subset L^q, H^{1-\alpha}\subset L^p$ in the second inequality and the interpolation inequality and Young's inequality in the last inequality.
Thus (b.1') is satisfied if $s_0\geq 2-\alpha$. By Lemmas 5.1 and 5.2 we have for $s=s_0,s_0-1$, $v\in H^{s_0+1}$
$$|\Lambda^{s}(u_v\cdot\nabla v)|=|\Lambda^{s}\nabla\cdot(u_vv)|\leq C|\Lambda^{s+1+\alpha}v||\Lambda^{1-\alpha}v|.$$Here $u_v=R^\bot v$.
Then (b.2) is satisfied for $s_0\geq1$. Now we  prove (b.3). Indeed for $v_1,v_2\in H^{2-\alpha}$ by Lemma 5.2
$$\aligned |\langle v_1-v_2,u_1\cdot\nabla v_1-u_2\cdot\nabla v_2\rangle|=&|\langle v_1-v_2,(u_1-u_2)\cdot\nabla v_1\rangle|
\\\leq& C|\Lambda^{2-\alpha}v_1||\Lambda^\alpha (v_1-v_2)||v_1-v_2|.\endaligned $$Here $u_i=R^\bot v_i$ and we used $div u_2=0$ in the first equality and Lemma 5.2 in the last inequality. Thus (b.3) is satisfied if $s_0\geq 2-\alpha$.
Now Theorems 3.2 and 4.1 imply the following results:
\vskip.10in
\th{Theorem 5.6} Fix $0<\alpha<1$. Assume that  $G$ satisfies (G.1)(G.2) with $s_0\geq 2-\alpha$. Assume that $\theta_0$ is an $H^{s_0}$-valued, $\mathcal{F}_0$-measurable random variable with $E|\Lambda^{s_0}\theta_0|^2<\infty$.

(i) Then local pathwise uniqueness holds and there exists a maximal strong solution $(\theta,(\tau_R)_{R\in\mathbb{N}},\zeta)$ of (5.2).

(ii) Moreover, for $G(\theta)f_i=\delta_{1i}\beta\theta$  there exists a positive deterministic funtion $\kappa:[1,\infty)\times (0,\infty)\rightarrow (0,\infty)$ defined by $\kappa(R,\beta^2):=(\frac{\beta^2}{4C_1})^{\delta(s_0)}\frac{1}{R}$, where $C_1,\delta(s_0)$ are as above, satisfying
$$\lim_{\beta^2\rightarrow\infty}\kappa(R,\beta^2)=\infty,$$
 such that whenever
$$|\Lambda^{s_0}\theta_0|^2\leq\kappa(R,\beta^2) \qquad P-a.s.,$$
then $$P(\zeta=\infty)\geq1-\frac{1}{R^{1/4}}$$
and $$P(\lim_{t\rightarrow\infty}|\Lambda^{s_0}\theta(t)|=0)\geq1-\frac{1}{R^{1/8}}.$$

\th{Remark 5.7}For $\alpha=1/2$, consider
$$d\theta+[A^\alpha\theta+u\cdot\nabla\theta]dt=\sum_{j=1}^mb_j\theta\circ dw_j(t),\eqno(5.3)$$
for $b_j\in\mathbb{R}$,  and independent $1$-dimensional Brownian motions $w_j$. Here $\circ$ means the Stratonoaich integral.
Consider the process
$$\beta(t)=e^{-\sum_{j=1}^m b_jw_j(t)}.$$
Then, the process $v(t)$ defined by the transformation
$$v(t)=\beta(t)\theta(t),$$
satisfies the equation with one coefficient depending on $\omega\in\Omega$
$$\frac{dv}{dt}+A^\alpha v+\beta^{-1}u_v\cdot\nabla v=0.\eqno(5.4)$$
Then by the same arguments as in the proof of  Theorem 3.1, we obtain a maximal strong solution $(\theta,(\tau_R)_{R\in\mathbb{N}},\zeta)$ to (5.3) starting from any point in $H^{\frac{3}{2}}$ and $\theta(\cdot\wedge\tau_R)\in C([0,\infty),H^{\frac{3}{2}})$. Then by It\^{o}'s formula, $v$ is a solution to (5.4) and $v(\cdot\wedge\tau_R)\in C([0,\infty),H^{\frac{3}{2}})$.
On the other hand, by the same arguments as in [CV06, Section 2], we obtain for fixed $\omega$ and any $T>0$, that there exists a constant $M=M(\omega,|\Lambda \theta_0(\omega)|)$ such that
$$\|v(t,\cdot)\|_\infty \leq M \textrm{ for } t\in [0,T].$$
Then $$\|\beta^{-1}u_v(t,\cdot)\|_{\rm{BMO}}\leq M_1(\omega,|\Lambda \theta_0|,T) \textrm{ for } t\in [0,T].$$
Hence by [KN09, Theorem 1.1], we obtain that there exists $\gamma(\omega,|\Lambda \theta_0|,T)>0$, such that
$$\|v(\cdot,t)\|_{C^\gamma(\mathbb{T}^2)}\leq C(\omega,|\Lambda \theta_0|,T).$$
Then by the same arguments as in the proof of [CV06, Theorem 10], we obtain for any $0<\beta<1$
$$\|v(\cdot,t)\|_{C^{1,\beta}(\mathbb{T}^2)}\leq C_1(\omega,|\Lambda \theta_0|,T)\textrm{ for } t\in [0,T].$$
By this a-priori bound and the local existence, we obtain
$$\zeta=\infty \quad P-a.s.,$$
which implies the existence and uniqueness of a global solution to (5.3).

\subsection{ Stochastic fractional Navier-Stokes equation in d-dimensions}

For $d\geq 2$ we consider the following $d$-dimensional stochastic fractional Navier-Stokes equation in a bounded open domain $O\subset \mathbb{R}^d$ with smooth boundary:
$$du+((u\cdot\nabla u)+(-\Delta)^\alpha u)dt=\nabla p+G(u) dW, \eqno(5.5)$$
$$\rm{div} u=0, u(0)=u_0,$$
$$u(t,x)=0, \quad (t,x)\in \mathbb{R}^+\times \partial O,$$
where $u=(u^1,...,u^d)$ represents the velocity field of the fluid, the pressure $p$ is an unknown scalar function.

When $\alpha=1,G=0$, (5.5) reduces to the usual Navier-Stokes equation. In the deterministic case for general $\alpha$ this equation has been studied by many authors (see [W05], [Z12] and the references therein). In [W05], the author  obtains global existence and uniqueness of solutions for small initial values and for $\alpha>\frac{1}{2}$. In [Z12], the author obtains the local existence and uniqueness of solutions by using a stochastic Lagrangian particle
trajectories approach for $\alpha=\frac{1}{2}$. Below we shall improve both results in an essential way as consequences of our main Theorems 3.2 and 4.1.

  Let $H:=\{f\in L^2(O;\mathbb{R}^d),\rm{div} f=0,\}$ and let $P$ be the orthogonal projection operator from $L^2(O)^d$ onto $H$. Then (5.5) can be rewritten as follows:
$$du+(P(u\cdot\nabla u)+(-\Delta)^\alpha u)dt=PG(u) dW, $$
$$\rm{div} u=0, u(0)=u_0,$$
$$u(t,x)=0, \quad (t,x)\in \mathbb{R}^+\times \partial O.$$
Set
$$Av:=-\Delta v\textrm{,}\quad v\in D(A):=H^2:=\{u\in H^{2,2}_0(O)^d:\rm{div} u=0\}.$$
 Then $A$ is an unbounded positive definite  self-adjoint operator and   $A^{-1}$ is compact on $H$. Then $H^{s}=\{u\in H^{s,2}_0(O)^d:\rm{div} u=0\}, s\in\mathbb{R}$. Set
 $$f(v):=P(v\cdot\nabla v)\textrm{,} \quad v\in H^{\frac{d}{2}+\delta_0},$$
for some $\delta_0>0$. Then we can easily find an $s_0\geq\frac{d}{2}+1-\alpha$ such that $f\in \Phi(H^{s_0};H)\cap \Phi(H^{s_0+1+\alpha};H^{s_0})$ and for every $\varphi\in \bigcap_{l=1}^\infty H^l$, $\langle f(\cdot),\varphi\rangle $ is continuous from $H^{s_0}$ to $\mathbb{R}$.

 Then, because $\langle u\cdot\nabla \Lambda^su,\Lambda^{s}u\rangle=0$ and by Lemma 5.3, we obtain that for any positive constant $a$ and $s\geq s_0, v\in H^{2s}$
$$\aligned |\langle av\cdot \nabla v,\Lambda^{2s}v\rangle|&=a|\langle \Lambda^{s}(v\cdot \nabla v)-v\cdot\nabla \Lambda^sv,\Lambda^{s}v\rangle|\\&\leq Ca|\Lambda^sv|\|\Lambda^sv\|_{L^{p_2}}\|\Lambda v\|_{L^{p_1}}\\&\leq Ca|\Lambda^sv||\Lambda^{s+\alpha}v||\Lambda^{1+\frac{d}{2}-\alpha} v|\leq \varepsilon|\Lambda^{s+\alpha}v|^2+C_1a^{2/\delta}|\Lambda^{1+d/2-\alpha}v|^{2+\frac{2}{\delta}},\endaligned$$
where $\frac{1}{p_1}+\frac{1}{p_2}=\frac{1}{2},  \delta:=\delta(s_0)=\frac{\alpha}{s+2\alpha-1-\frac{d}{2}}$ and we used $H^{d/2-\alpha}\subset L^{p_1}, H^\alpha\subset L^{p_2}$ in the second inequality and we used the interpolation inequality and Young's inequality in the last inequality. Thus (b.1') is satisfied if $s_0\geq 1+\frac{d}{2}-\alpha$. By Lemmas 5.1 and 5.2 we have for $s=s_0,s_0-1$, $v\in H^{s_0}$
$$|\Lambda^{s}(v\cdot\nabla v)|\leq C(|\Lambda^{s+1+\alpha}v||\Lambda^{\frac{d}{2}-\alpha}v|
+|\Lambda^{s+1+\alpha}v||\Lambda^{(\frac{d}{2}-\alpha)\vee 1}v|)\leq C|\Lambda^{s+1+\alpha}v||\Lambda^{(\frac{d}{2}-\alpha)\vee 1}v|.$$
Then (b.2) is satisfied for $s_0\geq1\vee(\frac{d}{2}-\alpha)$. Now we only need to prove (b.3). Indeed, for $v_1,v_2\in H^{\frac{d}{2}+1-\alpha}$,
$$\aligned |\langle v_1-v_2,v_1\cdot\nabla v_1-v_2\cdot\nabla v_2\rangle|=&|\langle v_1-v_2,(v_1-v_2)\cdot\nabla v_1\rangle|
\\\leq& C|\Lambda^{1+\frac{d}{2}-\alpha}v_1||\Lambda^\alpha (v_1-v_2)||v_1-v_2|,\endaligned $$
where we used $div v_2=0$ in the first equality and we used Lemma 5.2 in the last inequality. Thus (b.3) is satisfied if $s_0\geq 1+\frac{d}{2}-\alpha$.
Now Theorems 3.2  and 4.1 imply the following results:
\vskip.10in
\th{Theorem 5.8} Fix $0<\alpha\leq1$. Assume that  $G$ satisfies (G.1)(G.2) with $s_0\geq 1+\frac{d}{2}-\alpha$. Assume that $\theta_0$ is an $H^{s_0}$-valued, $\mathcal{F}_0$-measurable random variable with $E|\Lambda^{s_0}\theta_0|^2<\infty$.

(i) Then local pathwise uniqueness holds and there exists a maximal strong solution $(\theta,(\tau_R)_{R\in\mathbb{N}},\zeta)$ of (5.5).

(ii) Moreover, for $G(\theta)f_i=\delta_{1i}\beta\theta$  there exists a positive deterministic funtion $\kappa:[1,\infty)\times (0,\infty)\rightarrow (0,\infty)$ defined by $\kappa(R,\beta^2):=(\frac{\beta^2}{4C_1})^{\delta(s_0)}\frac{1}{R}$, where $C_1,\delta(s_0)$ are as above, satisfying
$$\lim_{\beta^2\rightarrow\infty}\kappa(R,\beta^2)=\infty,$$
such that whenever
$$|\Lambda^{s_0}\theta_0|^2\leq\kappa(R,\beta^2) \qquad P-a.s.,$$
then $$P(\zeta=\infty)\geq1-\frac{1}{R^{1/4}}$$
and $$P(\lim_{t\rightarrow\infty}|\Lambda^{s_0}\theta(t)|=0)\geq1-\frac{1}{R^{1/8}}.$$

\subsection{Fractional K-P-Z equation}
Consider the following equation on the $\mathbb{T}^d$:
$$du+(\lambda|\nabla u|^2+(-\Delta)^\alpha u)dt=G(u) dW, \eqno(5.6)$$
$$ u(0)=u_0,$$
where $\lambda$ is a constant.
Originally,
this equation was proposed as a model of surface
growth in [KPZ86]. However, it was later realized that it is a universal object that
describes the fluctuations of a number of strongly interacting models of statistical
mechanics. In [Ha13] the author introduces a new concept of solution to the KPZ equation when the stochastic perturbation is space time white noise. In this section we consider the equation driven by multiplicative trace-class noise.

Let $H:=\{ f\in L^2(\mathbb{T}^d),\int fd\xi=0\}$ and
$$Av:=-\Delta v\textrm{,}\quad v\in D(A):=H^2= H^{2,2}(\mathbb{T}^d).$$
 Then $A$ is an unbounded positive definite  self-adjoint operator and   $A^{-1}$ is compact on $H$.  Hence $H^{s}= H^{s,2}(\mathbb{T}^d), s\in\mathbb{R}$. Set
 $$f(v):=\lambda|\nabla v|^2\textrm{,} \quad v\in H^{1+\frac{d}{2}}.$$
Then we can easily find an $s_0\geq\frac{d}{2}+1$ such that $f\in \Phi(H^{s_0};H)\cap \Phi(H^{s_0+1+\alpha};H^{s_0})$ and for every $\varphi\in \bigcap_{l=1}^\infty H^l$, $\langle f(\cdot),\varphi\rangle $ is continuous from $H^{s_0}$ to $\mathbb{R}$.
We obtain that for any positive constant $a, \alpha>\frac{1}{2}$ and $s\geq s_0>1+\frac{d}{2}, v\in H^{2s}$
$$\aligned |\langle a |\nabla v|^2,\Lambda^{2s}v\rangle|\leq&a |\Lambda^{s-\alpha}|\nabla v|^2||\Lambda^{s+\alpha}v|
\\\leq&a C|\Lambda^{s-\alpha+1} v|\|\nabla v\|_{L^\infty}|\Lambda^{s+\alpha}v|\\\leq&C_1a^{\frac{2}{\delta}} |\Lambda^{s_0} v|^{\frac{2}{\delta}+2}+\varepsilon|\Lambda^{s+\alpha}v|^2,\endaligned$$
where $\delta:=\delta(s)=\frac{2\alpha-1}{s+\alpha-s_0}$ and we used Lemma 5.1 in the second inequality and $H^{s_0-1}\subset L^\infty$, the interpolation inequality and Young's inequality in the last inequality.
Thus (b.1') is satisfied if $s_0> 1+\frac{d}{2}$. By Lemmas 5.1 and 5.2 we have for $s=s_0,s_0-1$, $v\in H^{s_0+1}$
$$|\Lambda^{s}(\nabla v\cdot\nabla v)|\leq C |\Lambda^{s+1} v|\|\nabla v\|_{L^\infty}.$$
Then (b.2) is satisfied. Now we only need to prove (b.3). Indeed, for $v_1,v_2\in H^{s_0}, s_0>1+\frac{d}{2}$,
$$\aligned &|\langle v_1-v_2, |\nabla v_1|^2-|\nabla v_2|^2\rangle|\\\leq& |\langle v_1-v_2, \nabla (v_1-v_2)\cdot \nabla v_1+\nabla v_2\cdot\nabla (v_1-v_2)\rangle|\\\leq&C(|\Lambda^{\frac{d}{2}+1}v_1|+|\Lambda^{\frac{d}{2}+1}v_2|+\|\nabla v_1\|_{L^\infty}+\|\nabla v_2\|_{L^\infty})|\Lambda^{1-\alpha}(v_1-v_2)||\Lambda^\alpha (v_1-v_2)|
\\\leq&C(|\Lambda^{s_0} v_1|^{\frac{2}{1-r}}+|\Lambda^{s_0} v_2|^{\frac{2}{1-r}})|v_1-v_2|^2+\varepsilon|\Lambda^\alpha (v_1-v_2)|^2,\endaligned$$
where $r=\frac{1-\alpha}{\alpha}$ and we used Lemmas 5.1 and 5.2 in the second inequality and $H^{s_0-1}\subset L^\infty$ and Young's inequality in the last inequality. Thus (b.3) is satisfied if $s_0>1+\frac{d}{2}$.
Now Theorems 3.2 and 4.1 imply the following results:
\vskip.10in
\th{Theorem 5.9} Fix $\alpha>\frac{1}{2}$.  Assume that  $G$ satisfies (G.1)(G.2) with $s_0> 1+\frac{d}{2}$. Assume that $\theta_0$ is an $H^{s_0}$-valued, $\mathcal{F}_0$-measurable random variable with $E|\Lambda^{s_0}\theta_0|^2<\infty$.

(i) Then local pathwise uniqueness holds and there exists a maximal strong solution $(\theta,(\tau_R)_{R\in\mathbb{N}},\zeta)$ of (5.6).

(ii) Moreover, for $G(\theta)f_i=\delta_{1i}\beta\theta$  there exists a positive deterministic funtion $\kappa:[1,\infty)\times (0,\infty)\rightarrow (0,\infty)$ defined by $\kappa(R,\beta^2):=(\frac{\beta^2}{4C_1})^{\delta(s_0)}\frac{1}{R}$, where $C_1,\delta(s_0)$ are as above, satisfying
$$\lim_{\beta^2\rightarrow\infty}\kappa(R,\beta^2)=\infty,$$
such that whenever
$$|\Lambda^{s_0}\theta_0|^2\leq\kappa(R,\beta^2) \qquad P-a.s.,$$
then $$P(\zeta=\infty)\geq1-\frac{1}{R^{1/4}}$$
and $$P(\lim_{t\rightarrow\infty}|\Lambda^{s_0}\theta(t)|=0)\geq1-\frac{1}{R^{1/8}}.$$

\subsection{Surface growth PDE with random noise}

We consider a model which appears in the theory of growth of surfaces, which describes an amorphous material deposited on an initially flat surface in high vacuum. The corresponding SPDE is the following equation on the interval $[0,L]$:
$$d\theta(t)=[-(\partial_\xi^4)^\alpha\theta(t)-\partial_\xi^2\theta(t)+\partial_\xi^2(\partial_\xi\theta(t))^2]dt+G(\theta(t))dW(t), \theta(0)=\theta_0.\eqno(5.7)$$

It is known in the literature that for the case $\alpha=1$, the (1-dimension) surface growth model has
some similar features of difficulty as the 3D Navier-Stokes equation. In particular, the uniqueness of weak
solutions for this model is still an open problem in both the deterministic and stochastic case.
We should remark that for the space time white noise
case, the existence of a weak martingale solution was obtained by Bl\"{o}mker, Flandoli and
Romito in [BFR09] for this model, and the existence of a Markov selection and ergodicity properties
were also proved there.

Let $$H:=L^2([0,L]),$$
$$Av:=\partial_\xi^4 v\textrm{,} \quad v\in D(A):=H_0^{4,2}([0,L]).$$
 Then  $A$ is an unbounded positive definite  self-adjoint operator and   $A^{-1}$ is compact on $H$. Then $H^{s}=H^{2s,2}_0([0,L]), s\in\mathbb{R}$. Set
$$f(v):=\partial_\xi^2v-\partial_\xi^2(\partial_\xi v)^2\textrm{,}\quad v\in H^1.$$
Then we can easily check  that for $s_0\geq 3/2$ and $\alpha>3/4$, $f\in \Phi(H^{s_0};H)\cap \Phi(H^{s_0+1+\alpha};H^{s_0})$ and for every $\varphi\in \bigcap_{l=1}^\infty H^l$, $\langle f(\cdot),\varphi\rangle $ is continuous from $H^{s_0}\times H^{s_0}$ to $\mathbb{R}$.
We have the following estimate for $\alpha>\frac{3}{4}, s\geq s_0,v\in H^{2s}$:
$$\aligned &|\langle\partial_\xi^2v-a\partial_\xi^2(\partial_\xi v)^2,A^sv\rangle|\\\leq&|\Lambda^{s+\frac{1}{2}}v|^2+Ca|\Lambda^{s+1-\alpha}(\partial_\xi v)^2||\Lambda^{s+\alpha}v|\\\leq&|\Lambda^{s+\frac{1}{2}}v|^2+Ca|\Lambda^{s+\frac{3}{2}-\alpha} v|\|\partial_\xi v\|_{L^\infty}|\Lambda^{s+\alpha}v|\\\leq&C|\Lambda^{s+\alpha}v|^{2-2\delta_0}|\Lambda^{s_0}v|^{2\delta_0}+Ca|\Lambda^{s_0} v|^{1+\delta}|\Lambda^{s+\alpha}v|^{2-\delta}\\\leq&C_1(a^{\frac{2}{\delta}}|\Lambda^{s_0}v|^{2+\frac{2}{\delta}}+|\Lambda^{s_0}v|^2)+\varepsilon|\Lambda^{s+\alpha}v|^2,\endaligned$$
where $\delta_0=\frac{\alpha-1/2}{s+\alpha-s_0}, \delta:=\delta(s)=\frac{2\alpha-\frac{3}{2}}{s+\alpha-s_0}$ and we used Lemma 5.1 in the second inequality and $H^{s_0}\subset H^{1,\infty}_0([0,L])$ and the interpolation inequality in the third inequality and Young's inequality in the last inequality. Thus (b.1') is satisfied. For $v\in H^{s_0+\frac{3}{2}}$ we have for $s=s_0,s_0-1$

$$|\Lambda^{s}f(v)|\leq |\Lambda^{s+1}v|+|\Lambda^{s+1}(\partial_\xi v)^2|\leq |\Lambda^{s+1}v|+C|\Lambda^{s+1+\frac{1}{2}} v|\|\partial_\xi v\|_{L^\infty}.$$
Then (b.2) is satisfied for $s_0\geq 1$ with $m_1=1$. Now we verify (b.3): for $v\in H^{s_0}, s_0\geq 1$,
$$\aligned&|\langle f(v_1)-f(v_2),v_1-v_2\rangle|\\\leq& |\Lambda^{\frac{1}{2}}(v_1-v_2)|^2+|\langle(\partial_\xi v_1)^2-(\partial_\xi v_2)^2,\partial_\xi^2(v_1-v_2)\rangle|\\\leq&|\Lambda^{\frac{1}{2}}(v_1-v_2)|^2+C(|\Lambda v_1|+|\Lambda v_2|+\|\partial_\xi v_1\|_{L^\infty}+\|\partial_\xi v_2\|_{L^\infty})|\Lambda^{1-\alpha+\frac{1}{2}} (v_1- v_2)||\Lambda^{\alpha}(v_1-v_2)|\\\leq&C(|\Lambda^{s_0} v_1|^{\frac{2}{1-r}}+|\Lambda^{s_0} v_2|^{\frac{2}{1-r}}+1)|v_1- v_2|^2+\varepsilon|\Lambda^\alpha(v_1-v_2)|^2,\endaligned$$
where $r=\frac{\frac{3}{2}-\alpha}{\alpha}$ and we used the interpolation inequality and Young's inequality in the last inequality.
\vskip.10in
Now Theorems 3.2 and  4.1 imply the following results:
\vskip.10in
\th{Theorem 5.10} Fix $\alpha>\frac{3}{4}$. Assume that  $G$ satisfies (G.1)(G.2) with $s_0\geq3/2$. Assume that $\theta_0$ is an $H^{s_0}$-valued, $\mathcal{F}_0$-measurable random variable with $E|\Lambda^{s_0}\theta_0|^2<\infty$.

(i) Then local pathwise uniqueness holds and there exists a maximal strong solution $(\theta,(\tau_R)_{R\in\mathbb{N}},\zeta)$ of (5.7).

(ii) Moreover, for $G(\theta)f_i=\delta_{1i}\beta\theta$  there exists a positive deterministic funtion $\kappa:[1,\infty)\times (0,\infty)\rightarrow (4C_1,\infty)$ defined by $\kappa(R,\beta^2):=(\frac{\beta^2}{4C_1}-1)^{\delta(s_0)}\frac{1}{R}$, where $C_1,\delta(s_0)$ are as above, satisfying
$$\lim_{\beta^2\rightarrow\infty}\kappa(R,\beta^2)=\infty,$$
 such that whenever for $\beta^2>4C_1$
$$|\Lambda^{s_0}\theta_0|^2\leq\kappa(R,\beta^2) \qquad P-a.s.,$$
then $$P(\zeta=\infty)\geq1-\frac{1}{R^{1/4}}$$
and $$P(\lim_{t\rightarrow\infty}|\Lambda^{s_0}\theta(t)|=0)\geq1-\frac{1}{R^{1/8}}.$$

\subsection{Stochastic reaction-diffusion equations}

Let $O$ be an open bounded domain in $\mathbb{R}^d$ with smooth boundary. Consider the following semilinear stochastic equation
$$d\theta+(-\Delta)^\alpha\theta+p(\theta)\theta=G(\theta)dW,\eqno(5.8)$$
where $p$ is a polynomial of degree $k$.

Let $$H:=L^2(O),$$
$$Av:=-\Delta v\textrm{,} \quad v\in D(A):=W_0^{2,2}(O).$$
Then $A$ is an unbounded positive definite  self-adjoint operator and   $A^{-1}$ is compact on $H$. Hence $H^{s}=W^{s,2}_0(O), s\in\mathbb{R}$. Set
$$f(v):=p(v)v\textrm{,}\quad v\in H^{s_0},$$
for $s_0>\frac{d}{2}$.
Then we can easily check  that for $s_0>\frac{d}{2}$, $f\in \Phi(H^{s_0};H)\cap \Phi( H^{s_0+1+\alpha};H^{s_0})$ and for every $\varphi\in \bigcap_{l=1}^\infty H^l$, $\langle f(\cdot),\varphi\rangle $ is continuous from $H^{s_0}$ to $\mathbb{R}$.

Since for any $m\in\mathbb{N}$, $v\in H^{s_0}\cap H^s, l\leq s,$ by Lemma 5.1 we have
$$\aligned|\Lambda^{l}(v^m)|\leq& C[|\Lambda^{l}v|\|v^{m-1}\|_{L^\infty}+|\Lambda^{l}(v^{m-1})|\|v\|_{L^\infty}]
\\\leq&C[|\Lambda^{l}v|\|v\|_{L^\infty}^{m-1}+|\Lambda^{l}(v^{m-1})|\|v\|_{L^\infty}],\endaligned\eqno(5.9)$$
we obtain the following estimate: for $s\geq s_0,v\in H^{2s}$
$$\aligned|\langle p(av)v, \Lambda^{2s} v\rangle|\leq& |\Lambda^{s-\alpha}(p(av)v)||\Lambda^{s+\alpha}v|\\\leq&C|\Lambda^{s-\alpha}v|[\|av\|_{L^\infty}^{k}
+1]|\Lambda^{s+\alpha}v|\\\leq &C[\|av\|_{L^\infty}^{k}
+1]|\Lambda^{s_0}v|^{1-r}|\Lambda^{s+\alpha}v|^{1+r}\\ \leq &C_1[\|av\|_{L^\infty}^{\frac{1}{\delta}}
+1]|\Lambda^{s_0}v|^{2}+\varepsilon|\Lambda^{s+\alpha}v|^2 ,\endaligned$$
where $r=(\frac{s-\alpha-s_0}{s+\alpha-s_0})\vee0, \delta:=\delta(s)=\frac{1-r}{2k}$ and we used (5.9) in the second inequality and the interpolation inequality and Young's inequality in the last two inequalities. Thus (b.1') is satisfied for $s_0>\frac{d}{2}$. Now for $s=s_0, s_0-1$, $v\in H^{s_0+1}$ we have
$$|\Lambda^{s}(p(v)v)|\leq C|\Lambda^{s}v|[\|v\|^{k}_{L^\infty}+1].$$Thus (b.2) is satisfied for $s_0>\frac{d}{2}$.
Now we verify (b.3) for $v\in H^{s_0}$:
$$\aligned|\langle b(v_1,v_1)-b(v_2,v_2),v_1-v_2\rangle|=&|\langle p(v_1)v_1-p(v_2)v_2,v_1-v_2\rangle|\\\leq &C(1+\|v_1\|_{L^\infty}^{k}+\|v_2\|_{L^\infty}^{k})|v_1-v_2|^2,\endaligned$$
where we used (5.9) in the last inequality.
Now Theorems 3.2 and 4.1 imply the following results:
\vskip.10in
\th{Theorem 5.11} Fix $0<\alpha\leq1$. Assume that  $G$ satisfies (G.1)(G.2) with $s_0> \frac{d}{2}$. Assume that $\theta_0$ is an $H^{s_0}$-valued, $\mathcal{F}_0$-measurable random variable with $E|\Lambda^{s_0}\theta_0|^2<\infty$.

(i) Then local pathwise uniqueness holds and there exists a maximal strong solution $(\theta,(\tau_R)_{R\in\mathbb{N}},\zeta)$ of (5.8).

(ii) Moreover, for $G(\theta)f_i=\delta_{1i}\beta\theta$  there exists a positive deterministic funtion $\kappa:[1,\infty)\times (0,\infty)\rightarrow (0,\infty)$ defined by $\kappa(R,\beta^2):=(\frac{\beta^2}{4C_1}-1)^{2\delta(s_0)}\frac{1}{R}$, where $C_1, \delta(s_0)$ are as above, satisfying
$$\lim_{\beta^2\rightarrow\infty}\kappa(R,\beta^2)=\infty,$$
 such that whenever
$$|\Lambda^{s_0}\theta_0|^2\leq\kappa(R,\beta^2) \qquad P-a.s.,$$
then $$P(\zeta=\infty)\geq1-\frac{1}{R^{1/4}}$$
and $$P(\lim_{t\rightarrow\infty}|\Lambda^{s_0}\theta(t)|=0)\geq1-\frac{1}{R^{1/8}}.$$

\end{document}